\documentclass[12pt]{article}
\usepackage{amssymb,amsmath}

\newtheorem{lemma}{Lemma}[section]
\newtheorem{theorem}{Theorem}[section]

\newtheorem{corollary}{Corollary}[section]
\newtheorem{remark}{Remark}[section]
\begin{document}
\title{ Exponential elliptic boundary value problems on
a solid torus in the critical of supercritical case }
\author{{\LARGE Athanase Cotsiolis and Nikos Labropoulos}\\
Department of Mathematics, University of Patras,\\ Patras 26110, Greece\\
\small{e-mails: cotsioli@math.upatras.gr and nal@upatras.gr}}
\maketitle

{\noindent}{{\bf Abstract:} In this paper we investigate the
behavior and the existence of positive and non-radially symmetric
solutions to nonlinear exponential elliptic model problems defined
on a solid torus $\overline{T}$ of $\mathbb{R}^3$, when data are
invariant under the group $G=O(2)\times I \subset O(3)$.  The
model problems of interest are stated below:
\begin{flushleft}
$\begin{array}{ll} {\bf(P_1)} \hspace{0.5cm} & \displaystyle
 \Delta\upsilon+\gamma=f(x)e^\upsilon, \;\upsilon>0\quad \mathrm{on}
\quad T, \quad\upsilon |_{_{\partial T}}=0.\end{array}$
\end{flushleft}
$\,\,\,$ and
\begin{flushleft}
$\begin{array}{ll}\bf{(P_2)} \hspace{0.5cm} & \displaystyle
\Delta\upsilon+a+fe^\upsilon=0,
\;\upsilon>0\quad \mathrm{on}\quad T,\\[1.3ex]  &\displaystyle
\frac{\partial \upsilon}{\partial n}+b+ge^\upsilon=0\quad
\mathrm{on} \quad{\partial T}.\end{array}$
\end{flushleft}
We prove that exist solutions which are $G-$invariant and these
exhibit no radial symmetries. In order to solve the above problems
we need to find the best constants in the Sobolev inequalities in
the exceptional case.\\

{\noindent}{\textbf{Key words:}$\;$ Exponential problems $\cdot$
Solid
torus $\cdot$ Sobolev inequalities $\cdot$ Critical of supercritical exponent }}\\

{\noindent}{\textbf{Mathematics Subject Classification (2000)}}
35J66 $\cdot$ 46E35 $\cdot$ 35B33

\section{Introduction}
\setcounter{equation}{0} $\,\,\,\,\,\,$ In recent years,
significant progress has been made on the analysis of a number of
important features of nonlinear partial differential equations of
elliptic and parabolic type.  The study of these equations has
received considerable attention, because of their special
mathematical interest and because of practical applications of the
torus in scientific research today.  For example in Astronomy,
investigators study the torus which is a significant topological
feature surrounding many stars and black holes \cite{Dav-Smi}. In
Physics the torus is being explored at the National Spherical
Torus Experiment (NSTX) at Princeton Plasma Physics Laboratory to
test the fusion physics principles for the spherical torus concept
at the MA level \cite{Kay}. In Biology some investigators
interested in circular DNA molecules detected in a large number of
viruses, bacteria, and higher organisms.  In this topologically
very interesting type of molecule, superhelical turns are formed
as the Watson-Crick double helix winds in a torus formation
\cite{Cox-Bat}.

Let the solid torus be represented by the equation
$$
\overline{T}=\left\{(x,y,z)\in \mathbb{R}^3: (
\sqrt{x^2+y^2}-l)^2+z^2\leq r^2, \, \,l>r>0\right\},
$$
and the subgroup $G=O(2)\times I$ of $\;O(3)$. Note that the solid
torus $\overline{T}\subset \mathbb{R}^3 $ is
invariant under the group $G$.\\
We consider the following nonlinear exponential elliptic boundary
problems
\begin{flushleft}
$\begin{array}{ll} {\bf(P_1)} \hspace{0.5cm} & \displaystyle
 \Delta\upsilon+\gamma=f(x)e^\upsilon, \;\upsilon>0\quad \mathrm{on}
\quad T, \quad\upsilon |_{_{\partial T}}=0\end{array}$
\end{flushleft}
$\,\,\,$ and
\begin{flushleft}
$\begin{array}{ll}\bf{(P_2)} \hspace{0.5cm} & \displaystyle
\Delta\upsilon+a+fe^\upsilon=0,
\;\upsilon>0\quad \mathrm{on}\quad T,\\[1.3ex]  &\displaystyle
\frac{\partial \upsilon}{\partial n}+b+ge^\upsilon=0\quad
\mathrm{on} \quad{\partial T},\end{array}$
\end{flushleft}
where $\Delta \upsilon =-\nabla^i\nabla_i \upsilon$ is the
Laplacian of $\upsilon$, $\frac{\partial}{\partial n}$ is the
outer unit normal derivative, $f, g$ are two smooth $G-$invariant
functions and $\gamma, a, b
\in \mathbb{R}$.\\
Clearly, a radially symmetric solution is a $G-$invariant
solution, for any subgroup $G$ of $O(n)$. The converse problem is
considered in this paper, that is we prove that there exist
positive solutions which are $G-$invariant  and non$-$radially
symmetrical if $G=O(2)\times I$.

Problems ${\bf (P_1)}$ and ${\bf (P_2)}$ own their origin to the
"Nirenberg Problem" posed in $1969-70$ in the following way:\\
Given a (positive) smooth function $f$ on $(\mathbb{S}^2, g_0)$
(close to the constant function, if we want), is it the scalar
curvature of a metric $g$ conformal to $g_0$? ($g_0$ is the
standard metric whose sectional curvature is $1$) (see \cite{Aub1}). \\
Recall that, if we write $g$ in the form $g=e^ug_0$, the problem
is equivalent to solving the equation:
$$
\Delta u + 2 = fe^u.
$$
Nirenberg Problem has been studied extensively and is completely
solved (see \cite{Aub}, \cite{Ber}, \cite{Mos1}, \cite{Kaz-War},
\cite{Che4}). Further, we refer the reader to \cite{Che-Din2,
Che-Din1}, \cite{Che}, \cite{Bre-Mer}, \cite{Cha-Gur-Yan},
\cite{Li1, Li3, Li2}, \cite{Bia1, Bia2}, \cite{Han-Li}, \cite{Li},
\cite{Aub-Cot2}, \cite{Fag1}, \cite{Adi-Gia}, \cite{Che-Li},
\cite{Mei}, in which the authors study this problem or its
generalization.

Best constants in Sobolev inequalities are fundamental in the
study of non-linear PDEs on manifolds, because of their strong
connection with the existence and the multiplicity of the
solutions of the corresponding problems (see for example
\cite{Oba}, \cite{Aub5}, \cite{Sch}, \cite{Heb-Vau}, \cite{Heb},
\cite{Dru1}, \cite{Aub-Cot1, Aub-Cot2}, \cite{Cot-Lab1, Cot-Lab2,
Cot-Lab3} and the references therein). It is also well-known, that
Sobolev embeddings can be improved in the presence of symmetries
in the sense that we obtain continuous embeddings in higher $L^p$
spaces, that it, allow us to solve equations with higher critical
exponents (see for example \cite{Lio}, \cite{ Din}, \cite{Cot},
\cite{Cot-Ili}, \cite{Heb-Vau}, \cite{Dem-Naz} \cite{Aub-Cot1,
Aub-Cot2}, \cite{Fag1,Fag2}, \cite{Cot-Lab1, Cot-Lab2, Cot-Lab3}
and the references therein). Especially, in our case we solve
problems with the highest supercritical exponent (critical of
supercritical).

Let:
$$
C_{0,G}^\infty =\left\{\upsilon\in C_0^\infty(T):\upsilon\circ
\tau=\upsilon\, ,\, \forall \: \tau \in G\right\},
$$

$$
C_{G}^\infty =\left\{\upsilon\in C^\infty(T):\upsilon\circ
\tau=\upsilon\, ,\, \forall \: \tau \in G\right\}
$$
and
$$
L_{G}^p =\left\{\upsilon\in L^p(T):\upsilon\circ \tau=\upsilon\,
,\, \forall \: \tau \in G\right\}
$$
that is, the spaces of  all $G-$invariant functions  under the
action of the group $G=O(2)\times I$.

We define the Sobolev space $H^p_{1, G}(T)$, $p\geq 1$  as the
completion of $C^\infty_G(T)$ with respect to the  norm
$$
\|\upsilon\|_{H^p_1}=\|\nabla \upsilon\|_p+ \|\upsilon\|_p
$$
and $\stackrel{\circ}{H}\!_{1, G}^p(T)$ as the closure of
$C^\infty_{0,G}(T)$ in $H^p_{1, G}(T)$.

In \cite{Cot-Lab1} we proved that for any $p\in[1,2)$ real, the
embedding $H^p_{1,G}(T)\hookrightarrow L^q_G(T)$ is compact for
$1\leq q<\frac{2p}{2-p}$, while the embedding
$H^p_{1,G}(T)\hookrightarrow L^\frac{2p}{2-p}_{G}(T)$ is only
continuous. Also, in \cite{Cot-Lab2} we proved that for any
$p\in[1,2)$ real, the embedding $H^p_{1,G}(T)\hookrightarrow
L^q_G(\partial T)$ is compact for $1\leq q<\frac{p}{2-p}$, while
the trace embedding $H^p_{1,G}(T)\hookrightarrow
L^\frac{p}{2-p}_G(\partial T)$ is only continuous. Additionally,
we observe that if $\frac{3}{2}<p<2$ then
$q=\frac{2p}{2-p}>6=\frac{2\cdot 3}{3-2}$ and $\tilde q
>\frac{p}{2-p}>4=\frac{2(3-1)}{3-2}$, that is the
exponents $q$ and $\tilde q$ are supercritical.

In this paper, we study the exceptional case when $p=n-k=3-1=2$.
In this case $H_{1,G}^{2}(T)\not \subset L_G^\infty (T)$, however,
when $ \upsilon \in H_{1, G}^{2} (T)$ we have $e^\upsilon \in
L^1_G(T)$, $e^\upsilon \in L_G^1(\partial T)$ and the exponent
$p=2$ is the critical of supercritical.

This paper is organized  as follows:\\
In Section 2, we recall some definitions and we  present   the two
lemmas on which are based the proofs of the theorems concerning
the best constants. Proofs of these lemmas are in Section 6.
Section 3 is devoted in  the presentation of results of the paper.
In Section 4, we determine the best constants $\mu$ and
$\tilde\mu$ of the inequalities:
              \[\int_{T} e^\upsilon dV\leq C exp \left[\mu \left\|\nabla
               \upsilon  \right\|^2_2+ \frac{1}{2\pi^2r^2l}\int_{T} \upsilon dV\right]\]
               and
\[\int_{\partial T} e^\upsilon dS\leq C
exp\left[\tilde\mu\left\|\nabla\upsilon\right\|^2_2+
\frac{1}{4\pi^2rl} \int_{\partial T} \upsilon dS\right],\] In
section 5, we use the above two inequalities, in order to solve
the nonlinear exponential elliptic problems ${\bf(P_1)}$ and
${\bf(P_2)}$. Concerning problem ${\bf(P_1)}$, we prove the
existence of solutions of the associated variational problem. We
study  problem ${\bf(P_2)}$ in the same way as the ${\bf(P_1)}$,
except its last part (case 4 of Theorem 3.4), which is based upon
the  method of upper solutions and lower solutions.
\section{Notations and Preliminary Results}
For completeness we cite some background
material and results from \cite{Cot-Lab2}.\\
 Let  $\mathcal{A}=\{(\Omega_i,\xi_i):i=1,2\}$ be an atlas on $T$
defined by
\begin{gather*}
\Omega_1=\{(x,y,z)\in T :(x,y,z)\notin H^+_{XZ}\},\\
\Omega_2=\{(x,y,z)\in T :(x,y,z)\notin H^-_{XZ}\}
\end{gather*}
where
\begin{gather*}
H_{XZ}^+=\{(x,y,z)\in \mathbb{R}^3 : x>0\, ,\, y=0 \}\\
H^-_{XZ}=\{(x,y,z)\in \mathbb{R}^3 : x<0\, ,\, y=0 \}
\end{gather*}
and $\xi_i:\Omega_i\to I_i\times D$,  $i=1,2$, with
$
I_1=(0,2\pi), \quad I_2=(-\pi,\pi),
$
$$
D=\{(t,s)\in \mathbb{R}^2 :t^2+s^2<1\}, \quad \partial
D=\{(t,s)\in \mathbb{R}^2 :t^2+s^2=1\},
$$
$\xi_i(x,y,z)=(\omega_i,t,s)$, $i=1,2$  with $
\cos\omega_i=\frac{x}{\sqrt{x^2+y^2}},\;
\sin\omega_i=\frac{y}{\sqrt{x^2+y^2}},$ where
$$
\omega_1=  \begin{cases}
\mathop{\rm arctan}\frac{y}{x},&x\neq 0 \\
 \pi/2,&x=0,y>0   \\
3\pi/2,&x=0\,,\,y<0
 \end{cases}\quad
\omega_2=\begin{cases}
\mathop{\rm arctan}\frac{y}{x}, &x\neq 0 \\
     \pi/2,&x=0,\;y>0 \\
    -\pi/2,&x=0\,,\,y<0
     \end{cases}
$$
and
$$
 t=\frac{\sqrt{x^2+y^2}-l}{r}\,,\quad  s=\frac{z}{r}\,,\quad  0\leq t,\,s \leq 1.
$$
The Euclidean metric $g$ on $(\Omega,\xi)\in \mathcal{A}$ can be
expressed as
$$
(\sqrt{g}\circ\xi^{-1})(\omega,t,s)=r^2(l+rt).
$$

For any $G-$invariant $\upsilon$ we
define the functions $\phi(t,s)=(\upsilon\circ \xi^{-1})(\omega,t,s)$.\\
Then we have:
\begin{equation}\label{eq2.1}
    \int_T e^\upsilon dV= 2\pi r^2 \int_D e^{\phi(t,s)}(l+rt)\,dt\,ds
 \;\;\;\;
 \end{equation}
 \begin{equation}\label{eq2.2}
 ||\nabla \upsilon||_{L^2(T)}^2=2\pi  \int_D |\nabla \phi(t,s)|^2
 (l+rt)\,dt\,ds
 \end{equation}
and
\begin{equation}\label{eq2.3}
 \int_{\partial T} e^\upsilon dS=2\pi r \int_{\partial D
 }e^{\phi(t,0)}(l+rt)\,dt,\;\;\;\;\;\;\;\;
 \end{equation}
where by $\phi$ we denote the extension of  $\phi$
on $\partial D$.\\
Consider a finite covering $(T_j)_{j=1,\dots,N}$, where
$$
 T_j=\{(x,y,z)\in \overline {T} : ( \sqrt{x^2+y^2}-l_j)^2+(z-z_j)^2<
 \delta^2_j\}
 $$
is a tubular neighborhood (an open small solid torus) of the orbit
$O_{P_j}$ of $P_j$ under the action of the group $G$.
$P_j(x_j,y_j,z_j) \in\overline{T} $ and $l_j=\sqrt{x^2_j+y^2_j}\,$
is the horizontal distance of the orbit $O_{P_j}$  from the axis
$z'z$ and $\delta_j=l_j\varepsilon_j$ for any $\varepsilon_j>0$.\\
Then the following lemmas hold:
%
\begin{lemma}\label{L2.1}
$\bold{1.}\,\,$ For all $\varepsilon >0$, there  exists a constant
$C_\varepsilon$, such that for all $\upsilon\in
\stackrel{\circ}{H}\!^2_{1,G}(T_j)$ the following holds:
$$
\int_{T_j}e^\upsilon dV\leq C_\varepsilon exp \left[(1+c
\varepsilon ) \frac{1}{32\pi^2 l_j}\left\|\nabla \upsilon
\right\|^2_2\right],
$$
where $c>0$.\\
$\bold{2.}\,\,$ For all $\varepsilon >0$, there exist constants
$C_\varepsilon$
    and $D_\varepsilon$, such that for all $\upsilon\in
    \mbox{$\stackrel{\circ}{H}\!^2_{1,G}$}$ the following holds:
    $$
    \int_T e^\upsilon dV \leq C_\varepsilon exp
    \left[\left(\frac{1}{32\pi^2(l-r)}+ \varepsilon \right)
    \left\|\nabla \upsilon \right\|^2_2+D_\varepsilon \left\|\upsilon
    \right\|^2_2\right].
    $$
In addition the constant $\frac{1}{32\pi^2(l-r)}$ is the best
constant for the above inequality.
\end{lemma}

\begin{lemma}\label{L2.2}
Let $\overline{T}$ be the solid torus, $2\pi^2r^2l$ be the volume
of $\,T$ and $4\pi^2rl$ \\be the volume of $\,\partial T$, then
for all $\varepsilon >0$ there
exists a constant $C_\varepsilon$ such that:\\
$\bold{1.}\,\,$For all functions $\upsilon\in \mathcal{H}_G$ the
following inequality holds
\begin{equation}
\label{eq3.6} \int_T e^\upsilon dV \leq C_\varepsilon exp
\left[(\mu+\varepsilon)\|\nabla\upsilon\|^2_2+
\frac{1}{2\pi^2r^2l}\int_T \upsilon dV\right]
\end{equation}
$\bold{2.}\,\,$ For all functions $\upsilon\in \mathcal{H}_G$ the
following inequality holds
\begin{equation}\label{eq3.7}
\int_{\partial T} e^\upsilon dS \leq C_\varepsilon exp
\left[(\mu+\varepsilon)\|\nabla\upsilon\|^2_2+
\frac{1}{4\pi^2rl}\int_ {\partial T }\upsilon dS\right],
\end{equation}
where, for the first inequality, $\mu=\frac{1}{32\pi^2(l-r)}$ if
$\mathcal{H}_G=\stackrel{\circ}{H}\!_{1,G}^2$ and\\
$\mu=\frac{1}{16\pi^2(l-r)}$ if $\mathcal{H}_G=H_{1,G}^2$.\\ For
the second inequality $\mu>\frac{1}{8\pi^2(l-r)}$ for all
$\upsilon\in H^2_{1,G}$.\\ The constant $\mu$ is the best constant
for the above inequalities.
\end{lemma}
\section{Statement of Results}

\subsection{Best constants on the solid torus}

We have the following theorem:
\begin{theorem}\label{T3.1}
Let $\overline{T}$ be the solid torus, $2\pi^2r^2l$ be the volume
of $\,T$ and \\$4\pi^2rl$ be the volume of $\,\partial T$, then
there exists a constant $C$ such that:\\
$\bold{1.}\,\,$For all functions $\upsilon\in \mathcal{H}_G$ the
following inequality holds
\begin{equation}\label{eq3.1}
\int_T e^\upsilon dV\leq C exp \left[\mu\|\nabla\upsilon\|^2_2+
\frac{1}{2\pi^2r^2l}\int_T \upsilon dV\right]
\end{equation}
$\bold{2.}\,\,$ For all functions $\upsilon\in \mathcal{H}_G$ the
following inequality holds
\begin{equation}\label{eq3.2}
\int_{\partial T} e^\upsilon dS \leq C exp
\left[\mu\|\nabla\upsilon\|^2_2+ \frac{1}{4\pi^2rl}\int_{\partial
T} \upsilon dS\right],
\end{equation}
where, for the first inequality, $\mu=\frac{1}{32\pi^2(l-r)}$ if
$\mathcal{H}_G=\stackrel{\circ}{H}\!_{1,G}^2$ and\\
$\mu=\frac{1}{16\pi^2(l-r)}$ if $\mathcal{H}_G=H_{1,G}^2$.\\ For
the second inequality $\mu>\frac{1}{8\pi^2(l-r)}$ for all
$\upsilon\in H^2_{1,G}$.\\ The constant $\mu$ is the best constant
for the above inequalities.
\end{theorem}
\begin{remark}
\emph{In \cite{Fag4} Faget proved that for a compact
$3-$dimensional manifold without boundary the first best constant
for inequality  (\ref{eq3.1}) is $\mu_3=\frac{2}{81\pi}$ and the
map: $H_1^3 \ni \upsilon \rightarrow e^\upsilon \in L^1$ is
compact.  Clearly, the best constant $\mu_3$ depends only on the
dimension 3 of the manifold. For the solid torus, we prove that
the first best constant for the same inequality (\ref{eq3.1})  is
$\mu=\frac{1}{32\pi^2(l-r)}$ and the map: $H_{1,G}^2 \ni \upsilon
\rightarrow e^\upsilon \in L_G^1$ is compact. In this case, the
best constant $\mu$ depends on the  geometry of the solid torus.}
\end{remark}
\begin{corollary}\label{C3.1}
For all  $\;\upsilon\in\stackrel{\circ}{H}\!_{1,G}^2\;$ such that
$\;\left\|\nabla \upsilon \right\|^2_2\leq2\pi(l+r)\;$  and for
all
\\$\alpha \leq 4\pi$ the following holds:
\begin{equation}\label{eq3.11}
\int_T e^{\alpha\upsilon^2} dV \leq C\;2\pi^2r^2l
\end{equation}
where the constant $C$ is independent of
$\upsilon\in\stackrel{\circ}{H}\!_{1,G}^2$. The constant
$\alpha\leq 4\pi$ is the best, in the sense that, if $\alpha
> 4\pi$ the integral in the inequality is finite but
it can be made arbitrary large by an appropriate choice of
$\upsilon$.
\end{corollary}

\begin{remark}
\emph{Corollary \ref{C3.1} is a special case of the result of
Moser \cite{Mos}}.
\end{remark}
\subsection{Resolutions of the  Problems }


For the problem
\begin{flushleft}
$\begin{array}{ll} {\bf(P_1)} \hspace{0.5cm} & \displaystyle
 \Delta\upsilon+\gamma=f(x)e^\upsilon, \;\upsilon>0\quad\mathrm{on}
\quad T, \quad\upsilon |_{_{\partial T}}=0.\end{array}$
\end{flushleft}
we have the theorem:
\begin{theorem}\label{T3.3}
Consider a solid torus $\overline T$ and the function $f$ continuous and $G-$invariant.\\
Then the problem ${\bf (P_1)}$ accepts a solution that belongs to
$C^\infty_G $, if one of the following holds:
 \begin{description}
    \item [] \emph{(a)} \hbox{$\;\;$$sup_T f<0\,$ if $\,\;\gamma<0$.}
    \item [] \emph{(b)} \hbox{$\;\;$$\int_T fdV<0\,$ and $\,sup_T f>0\,$ if$\,\;\gamma=0$.}
    \item [] \emph{(c)} \hbox{$\;\;$ $sup_T f>0\,$ if $\;\,0<\gamma<\frac{8(l-r)}{lr^2}$.}
 \end{description}
\end{theorem}
For the problem
\begin{flushleft}
$\begin{array}{ll}\bf{(P_2)} \hspace{0.5cm} & \displaystyle
\Delta\upsilon+a+fe^\upsilon=0,
\;\upsilon>0\quad \mathrm{on}\quad T,\\[1.3ex]  &\displaystyle
\frac{\partial \upsilon}{\partial n}+b+ge^\upsilon=0\quad
\mathrm{on} \quad{\partial T},\end{array}$
\end{flushleft}
we have the next theorem:
\begin{theorem}\label{T3.4}
Consider a solid torus $\overline T$ and the smooth functions $f$,
$g$  $G-$ invariant and not both identical $0$. If $a,b\in
\mathbb{R}$ and $R= 2\pi^2r^2l a+ 4\pi^2rl b$, the problem ${\bf
(P_2)}$ accepts a solution that belongs to $C^\infty_G $ in each
one of the following cases :\\
   $\mathbf{1.}$ If $a=b=0$ the necessary and sufficient condition
          is $f$ and $g$ not both $\geq 0$ and that $\int_T f
          dV+\int_{\partial T} g dS>0$.\\
   $\mathbf{2.}$ If $a\geq 0$ and $b\geq 0$, $f$,$\;g$ not both
          $\geq 0$ everywhere and $\;0<R < 4\pi ^2 (l - r)$. Particularly, if
          $\;g=0$ then we can substitute the last condition with $\;0<R <
          8\pi ^2 (l - r)$.\\
   $\mathbf{3.}$ If $\;R>0$ (respectively $\;R<0$ ) it is necessary
          that $\;f$,$\;g$ not both $\geq 0$ everywhere (respectively
          $\;\leq0$). Then there exists a solution of the problem in each
          one of the following cases:
                 \begin{description}
                          \item[] \emph{(a)}  \hbox{$\;a<0$,$\;b>0$,$\;f<0$,$\;g\leq0$ and
                                                    $b<\frac{l-r}{lr}$
                                                    if $\;g\not\equiv0$ or $b<\frac{2(l-r)}{lr}$ if
                                                    $\;g\equiv0$}.\\
                          \item[] \emph{(b)}  \hbox{$\;a>0$,$\;b<0$,$\;f\leq0$,$\;g<0$
                                                     and $\;a<\frac{2(l-r)}{lr^2}$}.\\
                          \item[] \emph{(c)}  \hbox{$\;a>0$,$\;b<0$,$\;f\geq0$,$\;g>0$ and
                                       $\;a<\frac{2(l-r)}{lr^2}$ if $\;g\not\equiv0\;$ or
                                       $\;a<\frac{4(l-r)}{lr^2}$
                                       if $g\equiv0\;$}.\\
                          \item[]\emph{(d)}  \hbox{$\;a<0$,$\;b>0$,$\;f>0$,$\;g\geq0$ and
                           $\;b<\frac{l-r}{lr}$.}
                 \end{description}
$\mathbf{4.}$ If $\;a\leq0$,$\;b\leq0$, not both $=0$, it is
necessary $\int_T f dV+\int_{\partial T} g dS>0$. Then there
exists a non empty subset $S_{f,\;g}$ of $\;\mathbb R^2_-=\{(a,
b)\neq (0, 0): a\leq 0,\;b\leq 0\}$ with the property that if
$(c,d)\in S_{f,\;g}$ then $(c',d\,')\in S_{f,\;g}$ for any $c'\geq
c$,$\;d'\geq d $ and  such that the problem $\bf{(P_2)}$ has a
solution if and only if $\;(a,b)\in S_{f,\;g}$.
$S_{f,\;g}$$=\mathbb R^2_-$ if and only if the functions $f$,$g$
are
 $\not\equiv 0$ and $\geq0$. For all $(a,b)\in \mathbb R^2_-$ there exist functions $f$ and $g$
such that $\int_T f dV+\int_{\partial T} g dS>0$ and $(a,b)\notin
S_{f,\;g}$.

\end{theorem}

\section{Proofs of the Theorem concerning the best constants }

{\bf{Proof of Theorem \ref{T3.1}}}.
$\bold{1.}\,\,$ We give a proof by contradiction based on Lemma \ref{L2.2}.\\
Assume that for any $C_\alpha  $, there exist $\upsilon _\alpha
\in \stackrel{\circ}{H}\!^2_{1,G}$ with $\int_T {\upsilon _\alpha
} dV = 0$ such that
\begin{equation}\label{eq2.6}
\int_T {e^{\upsilon _\alpha  } dV > } \,C_\alpha  \exp \left( {\mu
\int_T {\left| {\nabla \upsilon _\alpha  } \right|} ^2 dV}
\right).
\end{equation}
Set $\phi_{\alpha }(t,s)=(\upsilon _{\alpha }\circ
\xi^{-1})(\omega,t,s)$. By ( \ref{eq2.6}) because of (\ref{eq2.1})
and  ( \ref{eq2.2}) we obtain sequentially
\begin{eqnarray*}
 2\pi r^2
\int_D e^{\phi _\alpha} (l+rt)  dtds\, >  \,C_\alpha
  \exp \left( 2\pi \mu \int_D \,\left| \nabla \phi _\alpha  \right|^2
( l + rt ) dtds \right),
\end{eqnarray*}
\begin{eqnarray*}
 2\pi r^2(l+r)
\int_D e^{\phi _\alpha}   dtds\, >  \,C_\alpha
  \exp \left( 2\pi \mu ( l - r )\int_D \,\left| \nabla \phi _\alpha  \right|^2
 dtds \right),
\end{eqnarray*}
and since $\mu=\frac{1}{32\pi^2(l-r)}$ (see part 1 of Lemma
\ref{L2.2}) we have
\begin{eqnarray*}
\int_D e^{\phi _\alpha}   dtds\, >  \,\frac{C_\alpha}{ 2\pi
r^2(l+r)}
  \exp \left( \frac{ 1}{16\pi}\int_D \,\left| \nabla \phi _\alpha  \right|^2
 dtds \right).
\end{eqnarray*}
The last inequality means that for any  $c_\alpha $, there exists
$\phi _\alpha   \in \mathop H\limits^ \circ \,_1^2 \left( D
\right)$ with \\$\int_D {\phi _\alpha  } dtds = 0$, such that
\begin{eqnarray*}
\int_D e^{\phi _\alpha}   dtds\, >  \,c_\alpha
  \exp \left(  \frac{ 1}{16\pi}\int_D \,\left| \nabla \phi _\alpha  \right|^2
 dtds \right),
\end{eqnarray*}
which is a contradiction, (see Theorem 1 in \cite{Che1}).\\

$\bold{2.}\,\,$ The proof of this part is similar to the proof of
the first one. Let us sketch it. Assume that for any $\tilde
C_\alpha $, there exist $\upsilon _\alpha \in \mathop {H\,_{1,G}^2
} $ with $\int_{\partial T }{\upsilon _\alpha  } dS = 0$ such that
\begin{equation}\label{eq2.130}
\int_{\partial T} {e^{\upsilon _\alpha  } dS > } \,\tilde C_\alpha
\exp \left( {\mu \int_T {\left| {\nabla \upsilon _\alpha  }
\right|} ^2 dV} \right)
\end{equation}
and define the function $\phi_{\alpha }$ as in the first part.\\
By ( \ref{eq2.130}) because of (\ref{eq2.2}) and  ( \ref{eq2.3})
we take  the inequality
\begin{eqnarray*}
\int_{\partial D} {e^{\phi _\alpha  } dt > } \,\tilde c_\alpha
\exp \left( {\mu \int_D {\left| {\nabla \phi _\alpha } \right|} ^2
dtds} \right),
\end{eqnarray*}
where $\tilde c_\alpha=\frac{\tilde C_\alpha}{2\pi r(l+r)}$.\\
The last inequality  is false (see Theorem 3 in \cite{Che3}) and the theorem is proved.\mbox{ }\hfill $\Box$\\
 \\
{\bf{Proof of Corollary \ref{C3.1}}}. Given $\varepsilon
>0$, let $(T_j)_{j=1,...,N}$ be a finite covering of $\bar{T}$, where
$$
T_j = \left\{ {Q \in \mathbb{R}^3 :d(Q,O_{P_j } ) < \delta _j
,\quad\delta _j=l_j\varepsilon_j}\quad\mathrm{and}
\quad\varepsilon _j \leqslant \varepsilon\right\}.
$$
For $\bar{T}$ we build a $G-$invariant partition of unity
$\left(h_j\right)_{j=1,...,N}$ relative to the $T_j$'s. If we
denote $\Phi = \upsilon \circ \xi _j^{ - 1}$, $ \Phi \in
\stackrel{\circ}{H}\!_1^2 (D)$,  for all $ \upsilon \in
\stackrel{\circ}{H}\!_{1,G}^2 $, following the same argument as in
the  Lemma \ref{L2.1} we obtain
\begin{eqnarray*}
\int_T {e^{\alpha  \upsilon ^2 } dV }&=& \int_T {\left( {\sum_{j =
1}^N {h_j } }
\right)e^{\alpha \upsilon ^2 } dV} \\
& =& \sum_{j = 1}^N {\int_{T_j }
{h_j e^{\alpha \upsilon  ^2 } dV} }  \\
& =& \sum_{j = 1}^N \int_{I\times D}{\left( {h_j  \circ \xi _j^{ -
1} } \right)}{e^{\alpha  \upsilon
^2  \circ \xi _j^{ - 1} } (\sqrt g  \circ \xi _j^{ - 1} )d\omega dtds}  \\
& =& 2\pi \sum_{j = 1}^N {\int_D {\left( {h_j  \circ \xi _j^{ - 1}
} \right)}
 e^{\alpha  \upsilon ^2  \circ \xi _j^{ - 1} } \delta _j^2 (l_j  + \delta _j t)dtds}  \\
& = &\frac{1}{\pi }\sum_{j = 1}^N {\int_D {\left( {h_j \circ \xi
_j^{ - 1} }
 \right)e^{\alpha  \Phi  ^2 } 2\pi ^2 \delta _j^2 l_j (1 + \varepsilon _j t)dtds} }  \\
& \le & \left( {1 + \varepsilon } \right)\frac{1}{\pi }\sum_{j =
1}^N {\int_D {\left( {h_j  \circ \xi _j^{ - 1} } \right)e^{\alpha
 \Phi ^2 }
 Vol\left( {T_j } \right)dtds} }  \le  \\
& \le &\left( {1 + \varepsilon } \right)\frac{1}{\pi }Vol\left( T
\right)\sum_{j = 1}^N
{\int_D {\left( {h_j  \circ \xi _j^{ - 1} } \right)e^{\alpha  \Phi ^2 } dtds} }  \\
& =& \left( {1 + \varepsilon } \right)2\pi r^2 l\int_D {\sum_{j =
1}^N {\left( {h_j
 \circ \xi _j^{ - 1} } \right)e^{\alpha  \Phi ^2 } dtds} }  \\
& =& \left( {1 + \varepsilon } \right)2\pi r^2 l\int_D {e^{\alpha
 \Phi ^2 }
 dtds}  \\
\end{eqnarray*}
 or
\begin{equation}\label{eq5.33}
\int_T {e^{\alpha  \upsilon  ^2 } dV \leqslant }(1 + \varepsilon )
2\pi r^2 l\int_D {e^{\alpha  \Phi ^2 } dtds}.
\end{equation}
Because of $\left\| {\nabla \upsilon } \right\|_2^2  \leqslant
2\pi (l + r)$ and (\ref{eq2.2}) we obtain $\left\| {\nabla \Phi }
\right\|_2  \leqslant 1$ and according to Theorem $2.47$ of
\cite{Aub1} for all $\Phi  \in \stackrel{\circ}{H}\!_1^2 (D)$ with
$\left\| {\nabla \Phi } \right\|_2  \leqslant 1$ and for any
$\alpha  \leqslant 4\pi $ the following inequality holds
\begin{equation}\label{eq5.34}
\int_D {e^{\alpha  \Phi ^2 } dtds \leqslant C\pi},
\end{equation}
where the constant $C$ is the same for all open and bounded
subsets of $\mathbb{R}^2 $.\\
Thus, from inequalities (\ref{eq5.33}) and (\ref{eq5.34}) we
obtain
\begin{equation}\label{eq5.35}
\int_T {e^{\alpha  \upsilon ^2 }  dV \leqslant } (1 + \varepsilon
)C2\pi^2r^2l.
\end{equation}
Suppose now that inequality (\ref{eq5.35}) does not hold for
$\varepsilon=0$. That is, there exists $\upsilon\in
\stackrel{\circ}{H}\!_{1,G}^2 $ with $\left\| {\nabla \upsilon }
\right\|_2^2 \leqslant 2\pi (l + r)$ and $\theta >0$ such that the
following inequality holds
\begin{equation}\label{eq5.36}
\int_T {e^{\alpha  \upsilon ^2 } dV \geq (1+\theta) C2\pi^2r^2l}
\end{equation}
By (\ref{eq5.36}), and because of (\ref{eq5.33}) we obtain
\begin{equation}\label{eq5.37}
\left( {1 + \varepsilon } \right)2\pi r^2 l\int_D {e^{\alpha
 \Phi ^2 } dtds \geq } \,\left( {1 + \theta }
\right)\,C2\pi ^2 r^2 l.
\end{equation}
Since (\ref{eq5.37}) holds for any $ \varepsilon  > 0$ we can
choose  $ \varepsilon $ such that $ \varepsilon  < \theta $ and
(\ref{eq5.37}) yields
\[
2\pi r^2 l\int_D {e^{\alpha  \Phi ^2 } dtds \geq } \,\,\frac{{1 +
\theta }} {{1 + \varepsilon }}\,C2\pi ^2 r^2 l
\]
or
\begin{equation}\label{eq5.39}
\int_{D } {e^{\alpha  \Phi ^2 } dtds
> } \,C\pi.
\end{equation}
But according to Theorem 2.47 of \cite{Aub1} for all $ \Phi \in
\stackrel{\circ}{H}\!_1^2 (D )$ the following
\[
 \int_{D} {e^{\alpha  \Phi  ^2 } dtds \leq } \,C\pi
\]
holds. Thus (\ref{eq5.39}) is false and the corollary  is
proved.\mbox{ }\hfill $\Box$

\section{Proofs of the Theorems concerning the problems }
{\bf{Proof of Theorem \ref{T3.3}}}. We see that if $f$ is a
constant the problem can be solved immediately. If $f=0$ and
$\gamma=0$, solutions are all the
constants. If $\gamma f >0$ the constant $ln(\gamma/f)$ is the solution. \\
Consider the functional
\[
I(\upsilon ) = \int_T {\left| {\nabla \upsilon } \right|^2 dV +
2\gamma \int_T {\upsilon dV} } ,
\]
the set
$$
A = \left\{ \upsilon  \in H_{1,G}^2 :\int_T {fe^\upsilon dV} =
\gamma Vol(T)\right \}
$$
and denote
\[
 \nu  = \mathop {\inf }\limits_{\upsilon \in A}
I(\upsilon ).
\]
If $\gamma>0$, in order $A\neq \emptyset$, it is necessary $f$ to
be somewhere positive, if $\gamma<0$ it's necessary $f$ to be
somewhere negative, and if $\gamma=0$ it is necessary $f$ to
change sign. In the following we accept that $f$ satisfies the
above necessary condition and it is not a constant.

$\bold{(a)}\;\;$$\gamma<0$ and $f$ negative everywhere.\\
Combining  Jensen's inequality:
\[
\frac{1} {{Vol\left( T \right)}}\int_T {\upsilon dV} \leqslant \ln
\left( {\frac{1} {{Vol\left( T \right)}}\int_T {e^\upsilon dV} }
\right)
\]
along with the following inequality:
\[
\frac{1} {{Vol\left( T \right)}}\int_T {e^\upsilon  dV} \leqslant
\frac{1} {{Vol\left( T \right)\sup f}}\int_T {f(x)e^\upsilon dV} =
\frac{\gamma } {{\sup f}}
\]
we obtain
\begin{equation}\label{eq6.7}
\int_T {\upsilon dV}  \leqslant Vol\left( T \right)\ln \left(
{\frac{\gamma } {{\sup f}}} \right)
\end{equation}
and thus
\[
 I(\upsilon ) \geqslant 2\gamma Vol\left( T \right)\ln
\left( {\frac{\gamma } {{\sup f}}} \right).
\]
From the last inequality we conclude that $\nu$ is finite.\\
Let $ \{ \upsilon _i \}\in A$  be a minimizing sequence of $I$,
that is $I(\upsilon_i)\rightarrow\nu$. If we take $ I(\upsilon _i
)\leqslant 1 + \nu$ we obtain
\[
\int_T {\left| {\nabla \upsilon _i } \right|^2 dV + 2\gamma \int_T
{\upsilon _i dV} } \leqslant 1 + \nu
\]
thus
\[
1 + \nu  - 2\gamma \int_T {\upsilon _i dV}  \geqslant \int_T
{\left| {\nabla \upsilon _i } \right|^2 dV} \geqslant 0
\]
and
\begin{equation}\label{eq6.8}
\int_T {\upsilon _i dV}  \leq \frac{1 + \nu }{2\gamma}.
\end{equation}
By (\ref{eq6.7}) and (\ref{eq6.8}) we obtain $ \left| {\int_T
{\upsilon _i dV} } \right| < C$, where $C$ is a constant.\\
In addition, we have
\[
 \int_T {\left| {\nabla \upsilon _i } \right|^2 dV}
\leqslant 1 + \nu  - 2\gamma \int_T {\upsilon _i dV} \leqslant 1 +
\nu  - 2\gamma C.
\]
Thus $ \{ \upsilon _i \}$ is bounded in $ H_{1, G}^2 (T)$ and
there exists a subsequence of $\upsilon_i$, denoted again by
$\upsilon_i$ and a function $\bar{\upsilon}$ such that:\\[0.15cm]
\indent $(a)$ $ \{ \upsilon _i\} \rightharpoonup \bar \upsilon$ on
$H_{1, G}^2 (T)$, (by Banach's Theorem),\\[0.15cm]
\indent$(b)$ $ \{ \upsilon _i \}  \to \bar \upsilon$ on $ L^2_G
(T)$, (by Kondrakov's Theorem),\\[0.15cm]
\indent $(c)$ $ \{ \upsilon _i \}  \to \bar \upsilon$ a.e., (by
Proposition 3.43 of \cite{Aub1}) and
\\[0.15cm]
\indent$(d)$ $\{ e^{\upsilon_i }\}  \to e^{\bar \upsilon }$ on $
L^1_G(T)$, (by Theorem \ref{T3.1}).\\[0.15cm]
\noindent From $(c)$ arises that $\bar{\upsilon}$ is $G-$invariant
and so $\bar \upsilon  \in A$, thus $ I(\bar \upsilon ) \geqslant \nu$.\\
From $(d)$ we conclude that
$$
 \left\| {\bar \upsilon }
\right\|_{H_1^2 }  \leqslant \mathop {\lim}\limits_{i \to \infty }
\inf\left\| {\upsilon _i } \right\|_{H_1^2 }  = \nu
$$
and by definition of $\nu$  we obtain  $ I(\bar \upsilon ) = \nu$.\\
Using the variation method we can prove that $\bar{\upsilon}$ is a
week solution of the corresponding Euler equation and, by the
regularization Theorem of \cite{Tol} and Theorem $3.54$ of
\cite{Aub1}, we conclude  that $ \bar \upsilon  \in C_G^\infty$.

$\bold{(b)}\;\;$$\gamma=0$ and $f$ changes sign.\\
In this case we need the extra condition $ \int_T {f(x)dV < 0}$,
because if we multiply the equation of the problem by
$e^{-\upsilon}$ and integrate we obtain

\[
\int_T {f(x)dV = \gamma \int_T {e^{ - \upsilon } } } dV - \int_T
{e^{ - \upsilon } } \left| {\nabla \upsilon } \right|^2 dV,
\]
the second part of this equality is negative.\\
Since $\gamma=0$,
\[
 I(\upsilon ) = \int_T {\left| {\nabla
\upsilon } \right|^2 dV}
\]
and considering $\int_T {\upsilon dV}  = 0$, if we define the set
\[
\tilde A =\left \{ {\upsilon  \in H_{1,G}^2 :\int_T {\upsilon dV}
= 0,\int_T {fe^\upsilon  } dV = 0} \right\}
\]
we will have
\[
\nu  = \mathop {\inf }\limits_{\upsilon  \in \tilde{A}} I(\upsilon
)\geq 0.
\]
In the following we work in the same way as in $(a)$.\\
Thus,  there exists a minimizing subsequence of $\upsilon_i$,
denoted again by $\upsilon_i$ that's converge on a function
$\bar{\upsilon}\in \tilde A$.\\
If $\kappa$ and $\lambda$ are the Lagrange multipliers, the Euler
equation is
\[
\Delta \bar \upsilon  + \kappa  = \lambda f(x)e^{\bar \upsilon }.
\]
Intergrading  by parts, because of $ \int_T {fe^{\bar \upsilon } }
dV = 0$, we obtain $ \kappa  = 0$ and for the function
$\bar{\upsilon}$ holds
\begin{equation}\label{eq6.9}
 \Delta \bar \upsilon  = \lambda f(x)e^{\bar \upsilon }.
\end{equation}
By equation (\ref{eq6.9}) we obtain that $\bar{\upsilon}$ is not
constant, because of $\int_T {f(x)dV < 0}$, and so $\lambda\not=
0$. In addition, multiplying the same equation by $ e^{ - \bar
\upsilon } $ and integrating by parts we obtain $ \lambda
\int_T {f(x)e^{\bar \upsilon } } dV < 0$ and then $\lambda > 0$.\\
Finally, is easy to check that the solution of the equation is $
\bar \upsilon  - \ln \lambda$.

$\bold{(c)}\;\;$$\gamma>0$ and $f$ somewhere positive.\\
Consider the same variation problem as in case $(a)$  and suppose
that $f$ is somewhere positive, which is the necessary condition
to be $A\not= \emptyset$, since $ \sup _T f > 0$. We have
\begin{equation}\label{eq6.10}
\gamma Vol\left( V \right) = \int_T {fe^\upsilon  dV \leqslant
\sup f\int_T {e^\upsilon  dV} }.
\end{equation}
In addition by Theorem \ref{T3.1} we have
\begin{equation}\label{eq6.11}
\int_T {e^\upsilon  dV}  \leqslant C\exp \left\{ {(\mu  +
\varepsilon )\int_T {\left| {\nabla \upsilon } \right|} ^2 dV +
\frac{1} {{Vol\left( V \right)}}\int_T {\upsilon dV} } \right\}.
\end{equation}
From (\ref{eq6.10}) and (\ref{eq6.11}) we obtain
\[
\gamma Vol\left( V \right)  \leqslant C\sup f\exp \left\{ {(\mu  +
\varepsilon )\int_T {\left| {\nabla \upsilon } \right|} ^2 dV +
\frac{1} {Vol(T)}\int_T {\upsilon dV} } \right\},
\]
\begin{eqnarray*}
\frac{\gamma Vol\left( V \right)}{C\sup f}  \leqslant \exp \left\{
{(\mu + \varepsilon )\int_T {\left| {\nabla \upsilon } \right|} ^2
dV + \frac{1} {Vol(T)}\int_T {\upsilon dV} } \right\},
\end{eqnarray*}
\begin{eqnarray*}
\ln \left( {\frac{{\gamma Vol\left( T \right)}} {{C\sup f}}}
\right) & \leqslant &(\mu  + \varepsilon )\int_T {\left| {\nabla
\upsilon } \right|^2 dV }+ \frac{1} {Vol(T)} \int_T {\upsilon dV},
\end{eqnarray*}
\begin{eqnarray*}
2\gamma Vol\left( T \right)\ln \left( {\frac{{\gamma Vol\left( T
\right)}} {{C\sup f}}} \right)  \leqslant 2\gamma Vol\left( T
\right)(\mu  + \varepsilon )\int_T {\left| {\nabla \upsilon }
\right|^2 dV }+ 2\gamma \int_T {\upsilon dV},
\end{eqnarray*}
\begin{eqnarray*}
2\gamma Vol\left( T \right)\ln \left( {\frac{{\gamma Vol\left( T
\right)}} {{C\sup f}}} \right) \!\leqslant \! 2\gamma Vol\left( T
\right)(\mu  + \varepsilon )\!\int_T {\left| {\nabla \upsilon }
\right|^2 dV \!+ I(\upsilon )}\! - \! \int_T {\left| {\nabla
\upsilon } \right|^2 dV},
\end{eqnarray*}
\begin{eqnarray*}
I(\upsilon ) \geqslant 2\gamma Vol(T)\ln \left( {\frac{{\gamma
Vol(T)}} {{C\sup f}}} \right) + \left[ {1 - 2\gamma Vol(T)(\mu  +
\varepsilon )} \right]\int_T {\left| {\nabla \upsilon } \right|^2
dV}
\end{eqnarray*}
or
\begin{equation}\label{eq6.12}
I(\upsilon ) \geqslant \left[ {1 - 2\gamma Vol(T)(\mu  +
\varepsilon )} \right]\int_T {\left| {\nabla \upsilon } \right|^2
dV}+C',
\end{equation}
where $ \mu  = \frac{1} {{32\pi ^2 (l - r)}}$ and $C'=2\gamma
Vol(T)\ln \left( {\frac{{\gamma Vol(T)}} {{C\sup f}}} \right) $.\\
So, for $\gamma  < \frac{{8(l - r)}} {{l\,r^2}}$, we have that $
I(\upsilon )$ is bounded bellow. \\
Thus if ${\upsilon_i}\in A$ is a minimizing sequence of $I$, by
equation (\ref{eq6.12}) we obtain that $ \left\| {\nabla \upsilon
_i }\right\|_2^2 \leqslant C_1$, and by equations (\ref{eq6.10})
and (\ref{eq6.11}) that $ \int_T {\upsilon _i } dV \geqslant C_2$,
where $C_1$ and $C_2$ are constants. Since $ \nu = \mathop {\inf
}_{\upsilon  \in A} I(\upsilon )$ and $ \mathop {\lim }_{i \to
\infty } I(\upsilon _i ) = \nu$ we may assume that $I(\upsilon _i
) < \nu + 1$ and so $ \int_T {\upsilon _i } dV \leqslant C_3$,
where $C_3$ is a constant. Thus $ \{ \upsilon _i \}$ is bounded in
$ H_{1, G}^2 (T)$ and then the rest of the proof is
the same as in case $(a)$.\mbox{ }\hfill $\Box$\\
{\bf{Proof of Theorem \ref{T3.4}}}. Following \cite{Che4}, let
$\upsilon \in C^\infty_G  (\overline{T})$ be a solution of
$\;{\bf(P_2)}$. We observe that integration by parts yields
\begin{flushleft}
$\begin{array}{lllll}  \hspace{1.5cm} & \displaystyle
  \int_T {(\Delta \upsilon  + a + fe^\upsilon  )dV}  = 0, \\[5ex]  &\displaystyle
   - \int_{\partial T} {\frac{{\partial\upsilon }}
{{\partial n}}dS}  + \int_T {(a + fe^\upsilon  )dV}  = 0,\\[5ex]  &\displaystyle
  \int_{\partial T} {(b + ge^\upsilon  )dS}  + \int_T {(a + fe^\upsilon  )dV}  = 0,
  \\[5ex]  &\displaystyle
  a\int_T {dV + b\int_{\partial T} {dS + } \int_T {fe^\upsilon  } dV}
    + \int_{\partial T} g e^\upsilon  dS = 0, \\[5ex]  &\displaystyle
  aVol(T) + bVol(\partial T) + \int_T {fe^\upsilon  } dV +
  \int_{\partial T} g e^\upsilon  dS = 0 \end{array}$
\end{flushleft}
namely
\begin{equation}\label{eq6.13}
K(\upsilon ) = aVol(T) + bVol(\partial T) + \int_T {fe^\upsilon  }
dV + \int_{\partial T} g e^\upsilon  dS=0.
\end{equation}
Multiplying by $e^{ - \upsilon } $ and integrating by parts also
implies
\begin{flushleft}
$\begin{array}{lllll}  \hspace{1.2cm} & \displaystyle
  \int_T {(e^{ - \upsilon } \Delta \upsilon  + ae^{ - \upsilon }  + f)dV}  = 0 \\[5ex]  &\displaystyle
   - \int_{\partial T} {e^{ - \upsilon } \frac{{\partial\upsilon }}
{{\partial n}}dS}  + \int_T {(ae^{ - \upsilon }  + f)dV}  - \int_T
{e^{ - \upsilon } \left| {\nabla \upsilon } \right|^2 dV}  = 0, \\[5ex]  &\displaystyle
  \int_{\partial T} {(e^{ - \upsilon } b + g)dS}  + \int_T {(ae^{ - \upsilon }
   + f)dV}  - \int_T {e^{ - \upsilon } \left| {\nabla \upsilon } \right|^2 dV}  = 0
\end{array}$
\end{flushleft}
namely
\begin{equation}\label{eq6.14}
a\int_T {e^{ - \upsilon } dV + } b\int_{\partial T} {e^{ -
\upsilon } dS}  + \int_T {fdV}  + \int_{\partial T} {gdS} - \int_T
{e^{ - \upsilon } \left| {\nabla \upsilon } \right|^2 dV}  = 0.
\end{equation}
Moreover, if $\upsilon  \in H_{1,G}^2 (T)$, according to
\cite{Aub1}, \cite{Tru} and \cite{Che3} and because of theorem
\ref{T3.1}, for any $q \geqslant 1$, $\upsilon \in L^q_G (T)$,
$\upsilon \in L^q_G (\partial T)$ and $e^\upsilon   \in
L^q_G (T)$.\\
Set
\[
I(\upsilon ) = \frac{1} {2}\int_T {\left| {\nabla \upsilon }
\right|^2 dV}  + a\int_T {\upsilon dV + } b\int_{\partial T}
{\upsilon dS}
\]
 and
 \[
 A =\left \{ \upsilon\in H_{1, G}^2 :K(\upsilon ) = 0\right\}.
\]
Our aim is the minimization of $I(\upsilon)$ on $A$.\\

$\bf{1.}\,$  Case $a=b=0$, $\int_T f dV+\int_{\partial T} g
dS>0$ and $f$ and $g$ not both $\geq 0$.\\
Since $f$ and $g$ are not both identically $0$, the solutions of
equation $\;{\bf(P_2)}$ are not constant functions. Hence if
$\upsilon$ is a solution we have
\begin{equation}\label{eq6.15}
\int_T {e^{ - \upsilon } \left| {\nabla \upsilon } \right|^2 dV}
> 0
\end{equation}
and then by (\ref{eq6.14}) and (\ref{eq6.15}) yield
\[
\int_T fdV + \int_{\partial T} {gdS}  > 0.
\]
Since $a = b = 0$ and $K(\upsilon ) = \int_T f e^\upsilon  dV +
\int_{\partial T} {ge^\upsilon  dS} $ in order $A = \{ \upsilon
\in H_1^2 (T):K(\upsilon ) = 0\}\not = \emptyset $
it's necessary $f$ and $g$ not to be both $\geq 0$.\\
Inversely, if $f$ and $g$ are not both $\geq 0$, we will prove
that $A\not=\emptyset$.\\
Because of
$$
\int_T {fdV}  + \int_{\partial T} {gdS}  > 0,
$$
we have
$$
\left\{f(T) \cup g(\partial T)\right\} \cap (0, + \infty ) \ne
\emptyset\quad \mathrm{ and }\quad \left\{f(T) \cup g(\partial
T)\right\} \cap ( - \infty ,0) \ne\emptyset .
$$
Define a $C^\infty$ function $\eta :[0, + \infty ) \to [0,1]$ such
that $\eta=1$ in $[0,1/2]$, $\eta=0$ in $[1, +\infty)$ and examine
the following two cases:

$\mathbf{(i) }$ $f$ changes sign on $T$.\\
There are two tori $T_1$ and $T_2$ contained in $T$ such that
$f>0$  on $ T_1$ and $f<0$ on $T_2$. Let the points $\, P_i , \;i
= 1, 2 \,$ belong to the central orbits $ \,O_{P_i } , \;i = 1,
2\,$ of $\,T_i , \;i = 1, 2\,$, respectively and let
\[
T_1  =\left \{ {(x,y,z) \in T:( {\sqrt {x^2  + y^2 }  - l_{P_1 } }
)^2 + ( {z - z_{P_1 } })^2  < \delta ^2 }\right \}
\]
and
\[
T_2  =\left \{ {(x,y,z) \in T :( {\sqrt {x^2  + y^2 }  - l_{P_2 }
} )^2  + ( {z - z_{P_2 } })^2  < \delta ^2 }\right \},
\]
where $ l_{P_i }  = \sqrt {x_{P_i }^2  + y_{P_i }^2 } ,\;i = 1,
2\;$ the horizontal distance of the orbit $ O_{P_i } ,i = 1,2$
from the
axis $z'z$.\\
Set
$$
\alpha  = \int_{T\backslash (T_1 \cup T_2 )} {fdV} +
\int_{\partial T} {gdS},
$$
and suppose that $\alpha\geq 0 $. Then  $$ \alpha _0 = \alpha  +
\int_{T_1 } f dV > 0.
$$
Consider the continuous function
\[\sigma (t) = \int_{T_2 } {f(P)\exp } \left[ {t\eta \left(
{\frac{d\left( {P,O_{P_2 }} \right)}{\delta }} \right)}
\right]dV,\;t \in \mathbb{R},\] where $d$ is the Euclidean
distance in $\mathbb{R}^3
$. \\
Since $\mathop {\lim }\limits_{t \to  + \infty } \sigma (t) = -
\infty $ and $\mathop {\lim }\limits_{t \to  - \infty } \sigma (t)
= 0$, there exists $t_0  \in \mathbb{R}$ such that $\sigma
(t_0 ) =  - \alpha _0 $. \\
Hence if we define the function $\upsilon  \in C^\infty _G (T)$ as
\[\upsilon (P) = \left\{ {\begin{array}{*{20}c}
   {t_0 \eta \left( {d\left( {P,O_{P_2 } } \right)/\delta } \right),\,\,\,\,P \in T_2 }  \\
   {\,\,\,\,\,0,\,\,\,\,\,\,\,\,\,\,\,\,\,\,\,\,\,\,\,\,\,\,\,\,\,
   \,\,\,\,\,\,\,\,\,\,\,\,\,\,\,\,\,\,\,P \notin T_2 \,\,\,\,}  \\
 \end{array} } \right.\]
by definition of $\sigma$ we obtain
 \[
 \sigma (t_0 ) = \int_{T_2 }
{f(P)e^{\upsilon (P)} } dV\] and then \[\int_T {f(P)} e^{\upsilon
(P)} dV =  - \alpha _0
\]
From the last equality we have
\begin{flushleft}
$\begin{array}{llll}  \hspace{1.5cm} & \displaystyle
  \int_{T_2 } {fe^\upsilon  } dV + \alpha  + \int_{T_1 }
  {fdV}  = 0, \\[5ex]  &\displaystyle
  \int_{T_2 } {fe^\upsilon  } dV + \int_{T\backslash (T_1
   \cup T_2 )} {fdV}  + \int_{\partial T} {gdS}  + \int_{T_1 }
   {fdV}  = 0,\\[5ex]  &\displaystyle
  \int_{T_2 } {fe^\upsilon  } dV + \int_{T\backslash T_2 } {fdV}
    - \int_{T_1 } {fdV}  + \int_{\partial T} {gds}
    + \int_{T_1 } {fdV}  = 0, \\[5ex]  &\displaystyle
  \int_{T_2 } {fe^\upsilon  } dV + \int_{T\backslash T_2 } {fdV}
    + \int_{\partial T} {gdS}  = 0 \end{array}$
\end{flushleft}
and from this by definition of $\upsilon $ we obtain
\begin{flushleft}
$\begin{array}{llll}  \hspace{1.5cm} & \displaystyle
  \int_{T_2 } {fe^\upsilon  } dV + \int_{T\backslash T_2 } {fe^\upsilon  dV}
   + \int_{\partial T} {ge^\upsilon  dS}  = 0, \\[5ex]  &\displaystyle
  \int_T {fe^\upsilon  } dV + \int_{\partial T} {ge^\upsilon  dS}  = 0. \end{array}$
\end{flushleft}
This means that $\upsilon  \in {\rm A}$ and hence
$A\not=\emptyset$.

$\mathbf{(ii) }$ $f$ does not change sign on $T$.\\
If $f \equiv 0\;$ and $g\;$ changes sing, following arguments of
the previous case, we construct a function $\upsilon  \in
C^\infty_G (\overline{T})$  such that $\int_{\partial
T}{ge^\upsilon  dS =
0}\; $, hence $K(\upsilon)=0\;$ and $A\not =\emptyset$.\\
If $\;f\not  \equiv 0\;$, let us suppose that $\;f \geqslant 0\;$
and $\;K(\upsilon ) = 0\;$. Then there exist $\;P_1  \in T\;$ and
$\;P_2 \in
\partial T\;$  such that $\;f(P_1 ) > 0\;$ and $\;g(P_2 ) < 0\;$. \\
Consider the tori
 \[
T_1  =\left \{ {(x,y,z) \in T:( {\sqrt {x^2  + y^2 }  - l_{P_1 } }
)^2 + ( {z - z_{P_1 } })^2  < \delta ^2 }\right\}
\]
and
\[
T_2  = \left\{ {(x,y,z) \in T:( {\sqrt {x^2  + y^2 }  - l_{P_2 } }
)^2 + ( {z - z_{P_2 } } )^2  < \delta ^2 }\right\},
\]
where $\,\delta\,$ is small enough, such that $\,\overline T_1
\cap
\overline T_2 = \emptyset \,$ , $f>0\,$ a.e. in $\,T_1\,$ and  $\,g<0\,$ a.e. in $\,T_2  \cap \partial T\,$.\\
Set
\[
\beta  = \int_{T\backslash T_1 } f dV + \int_{\partial T\backslash
T_2 } g dS + \int_{\partial T \cap T_2 } {g(P)\exp \left[ {t\eta
\left( {\frac{2d\left( {P,O_{P_2 } } \right)}{\delta }} \right)}
\right]} \,dS
\]
and choose $\,t\,$ large enough such that $\,\beta  < 0\,$.\\
Denote
\[
T_{2(\delta /2)}  = \left\{ {(x,y,z) \in T:( {\sqrt {x^2  + y^2 }
- l_{P_2 } } )^2  + ( {z - z_{P_2 } } )^2  < \left( {\frac{\delta
}{2}}\right)^2 }\right \}
\]
and define a function $ \vartheta  \in C^\infty  (\overline T),
\,0 \leqslant \vartheta (P) \leqslant 1 $ such that $\vartheta =1$
in a neighborhood of $\partial T \cap T_{2(\delta /2)} $,
$\vartheta =0$ out of $T_2$ and its support $K$ to have small
enough measure such that the following holds
\[
\gamma  = \int_K {f(P)\exp \left[ {t\vartheta (P)\eta \left(
{\frac{2d\left( {P,O_{P_2 } } \right)}{\delta }} \right)} \right]}
dV - \int_K {f(P)} dV <  - \beta.
\]
Consider now the continuous function
\[
\hat \sigma (t) = \int_{T_1 } {f(P)\exp \left[ {t\eta \left(
{\frac{d\left( {P,O_{P_1 } } \right)}{\delta }} \right)} \right]}
dV,\; t \in \mathbb{R}.
\]
Since $f \geqslant 0, \;f\not  \equiv 0$, $\;\mathop {\lim
}\limits_{t \to  - \infty } \hat \sigma (t) = 0\;$ and $\;\mathop
{\lim }\limits_{t \to  + \infty } \hat \sigma (t) =  + \infty \;$
there exists $\;t' \in \mathbb{R}\;$ such that $\;\hat \sigma (t')
= - (\beta + \gamma )\;$, that is
\begin{equation}\label{eq6.16}
   \int_{T_1 } {f(P)\exp \left[ {t'\eta \left( {\frac{d\left( {P,O_{P_2 } }
   \right)}{\delta }} \right)} \right]} dV =  - (\beta  + \gamma ) >
   0.
\end{equation}
Define now the function $\upsilon  \in C^\infty  (\overline T)$ by
\[
\upsilon (P) = \left\{ {\begin{array}{*{20}c}
  \;\;\;\;\;\;\; {t'\eta \left( {d\left( {P,O_{P_1 } } \right)/\delta } \right),\,P \in T_1 }  \\[2ex]
   {t\vartheta (P)\eta \left( {2d\left( {P,O_{P_2 } } \right)/\delta } \right),\,\,P \in T_2 }  \\[2ex]
   {0\,\,\,\,\,\,\,\,\,\,\,\,\,\,\,\,\,\,\,\,\,\,\,\,
   \,\,\,\,\,\,\,\,\,\,\,\,\,\,\,\,\,\,\,,P \notin T_1  \cup T_2 }  \\
 \end{array} } \right.
\]
We have
\begin{eqnarray}\label{eq6.17}
  \beta  &=& \int_{T\backslash T_1 } f dV + \int_{\partial T\backslash T_2 } g dS
  + \int_{\partial T \cap T_2 } {ge^\upsilon  } \,dS \nonumber \\
  & =& \int_{T\backslash T_1 } f dV + \int_{\partial T\backslash T_2 }
   g e^\upsilon  dS + \int\limits_{\partial T \cap T_2 } {ge^\upsilon  } \,dS \nonumber \\
& =& \int_{T\backslash T_1 } f dV + \int_{\partial T} g e^\upsilon
dS
\end{eqnarray}
and
\begin{eqnarray}\label{eq6.18}
  \gamma  &=& \int_{T_2 } {f(P)\exp \left[ {t\vartheta (P)\eta \left( {\frac{2d\left( {P,O_{P_2 } }
  \right)}{\delta }} \right)} \right]} dV \nonumber \\
  && - \int_{T_2 \backslash K} {f(P)\exp \left[ {t\vartheta (P)\eta \left( {\frac{2d\left( {P,O_{P_2 } }
   \right)}{\delta} } \right)} \right]} dV - \int_K {f(P)} dV\nonumber \\
  & =& \int_{T_2 } f e^\upsilon  dV - \int_{T_2 \backslash K} f dV -
   \int_K f dV.
\end{eqnarray}
By (\ref{eq6.16}), (\ref{eq6.17}) and (\ref{eq6.18}) we now obtain
$$
\int_{T\backslash T_1 } f dV - \int_{T_2 \backslash K} f dV -
\int_K f dV + \int_{T_1 } f e^\upsilon dV + \int_{T_2 } f
e^\upsilon  dV + \int_{\partial T} g e^\upsilon  dS = 0,
$$
$$
\int_{T\backslash T_1 } f dV - \int_{T_2 } f dV +
\int_{T_1 } f e^\upsilon dV + \int_{T_2 } f e^\upsilon  dV +
\int_{\partial T} g e^\upsilon  dS = 0,
$$
$$
 \int_{T\backslash
(T_1 \cup T_2 )} f dV + \int_{T_1 \cup T} f e^\upsilon  dV +
\int_{\partial T} g e^\upsilon  dS = 0, $$ $$ \int_{T\backslash
(T_1 \cup T_2 )} f e^\upsilon  dV + \int_{T_1  \cup T_2 } f
e^\upsilon dV + \int_{\partial T} g e^\upsilon  dS = 0,
$$
$$
\int_T f e^\upsilon  dV + \int_{\partial T} g e^\upsilon  dS = 0.
$$
Hence $\upsilon  \in {\rm A}$ and $A\not=\emptyset$.\\
We observe that if $K(\upsilon ) = 0$ then $K(\upsilon  + c) = 0$
for any constat $c$. So we can suppose that $\int_T {\upsilon dV}
= 0$ for any $\upsilon  \in {\rm A}$.\\
Set
\[\mu  =
\mathop {\inf }\limits_{\upsilon  \in A}\left\{ {\int_T {| {\nabla
\upsilon } |^2 } dV:\int_T {\upsilon dV} = 0} \right\} \geqslant
0.
\]
Let $\{ \upsilon _i \}$ be a minimizing sequence. Since $\sup _i
(\left\| {\nabla \upsilon _i } \right\|^2 ) <  + \infty $, this is
bounded in $H_{1,G}^2 (T)$. Thus there exists a subsequence $\{
\upsilon _i \}$ and a function $\upsilon  \in H_{1,G}^2 (T)$ such
that:\\[0.15cm]
\indent $(a)$ $ \{ \upsilon _i\} \rightharpoonup  \upsilon$ on
$H_{1,G}^2 (T)$, (by Banach's Theorem),\\[0.15cm]
\indent $(b)$ $ \{ \upsilon _i \}  \to \upsilon$ on $ L^q_G (T),
q\geq 1$, (by Kondrakov's Theorem),\\[0.15cm]
\indent $(c)$ $ \{ \upsilon _i \}  \to  \upsilon$ a.e., (by
Proposition 3.43 of \cite{Aub1}),\\[0.15cm]
\indent $(d)$ $ \{e^{ \upsilon_i }\}  \to e^{ \upsilon }$ (by
Theorem \ref{T3.1}) and\\[0.15cm]
\indent $(e)$ $ \{ \upsilon _i \}  \to  \upsilon$ a.e., on
$\partial T$ and $\{ e^{ \upsilon_i } \} \to e^{ \upsilon }$ on
$L^q_G(\partial T)$,\\[0.15cm]
\noindent where by $ \upsilon_i$ and $\upsilon$ we denote the
trace of $ \upsilon_i$ and  $\upsilon$ on $\partial T$,
respectively (by Theorem 4 of \cite{Che3}).  \\
The latter implies
\[
 \mathop {\lim }\limits_{i \to \infty }
\int_T {\upsilon _i dV = } \int_T {\upsilon dV}
\]
and
\[
\mathop {\lim }_{i \to \infty } \left( {\int_T f e^{\upsilon _i }
dV + \int_{\partial T} g e^{\upsilon _i } dS} \right) = \int_T f
e^\upsilon  dV + \int_{\partial T} g e^\upsilon  dS = 0.
\]
From the last two equalities along with $(c)$ arises that
$\upsilon \in {\rm A}$ and $\int_T {\upsilon dV}  = 0$, hence, by
definition of $\mu$, $\left\| {\nabla \upsilon } \right\|_2^2
\geqslant \mu $.\\
From $(b)$ and using Theorem $3.17$ of \cite{Aub1} we obtain
\[
\left\| {\nabla \upsilon } \right\|_2^2 \leqslant \mathop
{\lim}\limits_{i \to \infty } \inf\left\| {\nabla \upsilon _i }
\right\|_2^2  = \mu.
\]
Thus, by definition of $\mu$, $\left\| {\nabla \upsilon }
\right\|_2^2 = \mu $ and the $ \inf\left\| {\nabla \upsilon _i }
\right\|_2^2$ is attained, where $\upsilon
_i \in {\rm A}$ and $\int_T {\upsilon _i dV = } 0$ .\\
If $\kappa$ and $\lambda$ are the Lagrange multipliers, the Euler
equation is
\begin{equation}\label{eq6.19}
\int_T {\nabla ^i \upsilon \nabla _i h} dV + \kappa \left( {\int_T
{fe^\upsilon  h} dV + \int_{\partial T} {ge^\upsilon  } hdS}
\right) + \lambda \int_T {hdV = 0},
\end{equation}
for all $ h \in H_1^2 (T)$. \\
Since $ K(\upsilon ) = 0$, for $h=1$ arises $\lambda=0$,  and for
$ h = \upsilon$, $\kappa \ne 0$. (If $\kappa  = 0$, $\left\|
{\nabla \upsilon } \right\|_2  = 0$ and since $\int_T {\upsilon dV
= 0} $, $\upsilon =0$ a.e. thus $ K(\upsilon ) > 0$, which is false).\\
According to Theorem $1$ of  \cite{Che3} the solution $ \upsilon
\in H_{1,G}^2 $  of (\ref{eq6.19}) is $C^\infty $ and if
satisfies:
\begin{equation}\label{eq6.20}
\left. \begin{gathered} \Delta \upsilon  + \kappa fe^\upsilon   =
0\quad\mathrm{in}\quad
{\rm T}\, \hfill \\
\frac{{\partial \upsilon }} {{\partial n}} + \kappa ge^\upsilon
= 0\quad\mathrm{on}\quad
\partial {\rm T} \hfill \\
\end{gathered}  \right\}
\end{equation}
Setting $h = e^{ - \upsilon } $ in (\ref{eq6.19}) we find
\begin{equation}\label{eq6.21}
\kappa  = \left( {\int_T f dV + \int_{\partial T} g dS} \right)^{
- 1} \int_T {\left| {\nabla \upsilon } \right|^2 e^{ - \upsilon }
} dV > 0
\end{equation}
and then $\upsilon  - \ln \kappa $ is a solution of
$\;{\bf(P_2)}$.\\

$\bf{2.}\,$  Case $a\geq 0$, $b\geq 0$, not both $\equiv 0$ and
$$
f^{ - 1} \Bigl( {\left( { - \infty ,0} \right)} \Bigr) \ne
\emptyset \quad \mathrm{or} \quad g^{ - 1} \Bigl( {\left( { -
\infty ,0} \right)} \Bigr ) \ne \emptyset .
$$
In this case we have
$$
R = aVol(T) + bVol(\partial T) > 0
$$
and by (\ref{eq6.13})
$$\int_T {fe^\upsilon  dV}  + \int_{\partial
T} {ge^\upsilon  dS}  < 0.
$$
Then, if $f$, $g$ are not
both $ \geqslant 0$, $A\not=\emptyset$ . \\
By Theorem \ref{T3.1} arises that, for all $\varepsilon>0$, there
exists a constant $ C_\varepsilon$ such that
\begin{equation}\label{eq6.22}
\int_T e^\upsilon dV \leq {C_\varepsilon  \exp \left[ {\left( {1 +
\varepsilon } \right)\frac{1} {{16\pi ^2 \left( {l - r}
\right)}}\left\| {\nabla \upsilon } \right\|_2^2  + \frac{1}
{{Vol(T)}}\int_T {\upsilon dV} } \right]}
\end{equation}
and
\begin{equation}\label{eq6.23}
\int_{\partial T} e^\upsilon dS \leq {  C_\varepsilon \exp \left[
{\left( {1 + \varepsilon } \right)\frac{1} {{8\pi ^2 \left( {l -
r} \right)}}\left\| {\nabla \upsilon } \right\|_2^2  + \frac{1}
{{Vol(\partial T)}}\int_{\partial T} {\upsilon dS} } \right]},
\end{equation}
for all $\upsilon\in H^2_{1,G}$.\\
From the  definitions of $K(\upsilon)$ and  $R$ and by
(\ref{eq6.13}) we obtain
\[
 R = \left| {\int_T {fe^\upsilon  } dV +
\int_{\partial T} g e^\upsilon  dS} \right| \leqslant \left(
{\mathop {\max }_{\overline T} \left| f \right|} \right)\int_T
{e^\upsilon } dV + \left( {\mathop {\max }_{\partial T} \left| g
\right|} \right)\int_{\partial T} {e^\upsilon  } dS
\]
and using (\ref{eq6.22}), (\ref{eq6.23}) we obtain
\begin{eqnarray*}
R &\leqslant& \left( {\mathop {\max }_{\overline T} \left| f
\right|} \right)C_\varepsilon  \exp \left[ {\left( {1 +
\varepsilon } \right)\frac{1} {{16\pi ^2 \left( {l - r}
\right)}}\left\| {\nabla
\upsilon } \right\|_2^2  + \frac{1}{{Vol(T)}}\int_T {\upsilon dV} } \right] \hfill \\
  && + \left( {\mathop {\max }_{\partial T}  \left| g \right|} \right)C_\varepsilon
\exp \left[ {\left( {1 + \varepsilon } \right)\frac{1} {{8\pi ^2
\left( {l - r} \right)}}\left\| {\nabla \upsilon } \right\|_2^2
+\frac{1}{{Vol(\partial T)}}\int_{\partial T} {\upsilon dS} } \right] \hfill \\
&\leqslant& \left( {\mathop {\max }\limits_{\overline T} \left| f
\right|} \right)C_\varepsilon
  \exp \left[ {\left( {1 + \varepsilon } \right)\frac{1}
{{8\pi ^2 \left( {l - r} \right)}}\left\| {\nabla \upsilon }
\right\|_2^2  + \frac{1}
{{Vol(T)}}\int_T {\upsilon dV} } \right] \hfill \\
  && + \left( {\mathop {\max }\limits_{\partial T}
\left| g \right|} \right)C_\varepsilon  \exp \left[ {\left( {1 +
\varepsilon } \right)\frac{1} {{8\pi ^2 \left( {l - r}
\right)}}\left\| {\nabla \upsilon } \right\|_2^2  + \frac{1}
{{Vol(\partial T)}}\int_{\partial T} {\upsilon dS} } \right]. \hfill \\
\end{eqnarray*}
The last inequality gives
\begin{equation}\label{eq6.24}
\mathop {\inf }_{\upsilon  \in A} \left\{ {\int_T {\upsilon dV + }
\left( {1 + \varepsilon } \right)\frac{1} {{8\pi ^2 \left( {l - r}
\right)}}Vol(T)\left\| {\nabla \upsilon } \right\|_2^2 } \right\}
= c_T (\varepsilon ) > - \infty
\end{equation}
and
\begin{equation}\label{eq6.25}
\mathop {\inf }\limits_{\upsilon  \in A} \left\{ {\int_{\partial
T} {\upsilon dS + } \left( {1 + \varepsilon } \right)\frac{1}
{{8\pi ^2 \left( {l - r} \right)}}Vol(\partial T)\left\| {\nabla
\upsilon } \right\|_2^2 } \right\} = c_{\partial T} (\varepsilon )
>  - \infty.
\end{equation}
By (\ref{eq6.24}), (\ref{eq6.25}) we obtain
\begin{eqnarray}\label{eq6.26}
  I(\upsilon ) &=& \frac{1}
{2}\int_T {\left| {\nabla \upsilon } \right|^2 dV}  + a\int_T
{\upsilon dV +
 } b\int_{\partial T} {\upsilon dS} \nonumber \hfill \\
& \geqslant& \left[ {\frac{1} {2} - \left( {1 + \varepsilon }
\right)\frac{1} {{8\pi ^2 \left( {l - r} \right)}}\left( {aVol(T)
+ bVol(\partial T)} \right)}
 \right]\left\| {\nabla \upsilon } \right\|_2^2\nonumber  \hfill \\
&& + ac_T (\varepsilon ) + bc_{\partial T} (\varepsilon ) \nonumber\hfill \\
& =& \left[ {\frac{1} {2} - \left( {1 + \varepsilon }
\right)\frac{1} {{8\pi ^2 \left( {l - r} \right)}}R}
\right]\left\| {\nabla \upsilon } \right\|_2^2  + ac_T
(\varepsilon ) + bc_{\partial T} (\varepsilon ).
\end{eqnarray}
If we assume $R < 4\pi ^2 (l - r)$ and if we choose $\varepsilon
> 0$ such that $ c = \frac{1} {2} - \left( {1 + \varepsilon }
\right)\frac{1} {{8\pi ^2 \left( {l - r} \right)}}R > 0 $, by
(\ref{eq6.26}) we conclude that $\mu  = \mathop {\inf
}\limits_{\upsilon \in A} I(\upsilon ) >  - \infty $. \\
Let $\{ \upsilon _i \} _{i \in \mathbb{N}} , \upsilon _i  \in A$
be a minimizing sequence of $I(\upsilon )$ such that
\begin{equation}\label{eq6.27}
\mu  \leqslant I(\upsilon _i ) \leqslant \mu  + 1
\end{equation}
for any $i \in \mathbb{N}$ (\ref{eq6.26}) and (\ref{eq6.27}) yield
\begin{eqnarray}\label{eq6.28}
  0 \leqslant \left\| {\nabla \upsilon _i } \right\|_2^2  &\leqslant& \frac{{I(\upsilon _i )
  - ac_T (\varepsilon ) - bc_{\partial T} (\varepsilon )}}
{c}\nonumber \leqslant \frac{{\mu  + 1 - ac_T (\varepsilon ) -
bc_{\partial T} (\varepsilon )}} {c} <  + \infty.
\end{eqnarray}
By (\ref{eq6.24}), (\ref{eq6.25}) and (\ref{eq6.26}) we also
obtain
 \begin{equation}\label{eq6.29}
\int_T {\upsilon _i dV}  \geqslant c_T (\varepsilon ) - \left( {1
+ \varepsilon } \right)\frac{1} {{8\pi ^2 \left( {l - r}
\right)}}Vol(T)\left( {\mu  + 1} \right) = C_T
\end{equation}
and
\begin{equation}\label{eq6.30}
\int_{\partial T} {\upsilon _i dS}  \geqslant c_{\partial T}
(\varepsilon ) - \left( {1 + \varepsilon } \right)\frac{1} {{8\pi
^2 \left( {l - r} \right)}}Vol(\partial T)\left( {\mu  + 1}
\right) = C_{\partial T}.
\end{equation}
By the definition of $I(\upsilon)$ and because of (\ref{eq6.27})
yields
\[a\int_T {\upsilon _i dV + } b\int_{\partial T}
{\upsilon _i dS}  \leqslant I(\upsilon _i ) \leqslant \mu  + 1.\]
The last relation, because of  (\ref{eq6.29}), (\ref{eq6.30})
gives us
\begin{equation}\label{eq6.31}
\int_T {\upsilon _i dV}  \leqslant \frac{\mu  + 1}{a }-
C_{\partial T}\quad \mathrm{if} \quad a\neq 0
\end{equation}
and
\begin{equation}\label{eq6.32}
  \int_{\partial T} {\upsilon _i dS}  \leqslant \frac{\mu  + 1}{b} - C_T\quad \mathrm{if} \quad b\neq
  0.
\end{equation}
By (\ref{eq6.27}), (\ref{eq6.28}), (\ref{eq6.31}) and
(\ref{eq6.32}) we have
\begin{equation}\label{eq6.33}
\left| {\int_T {\upsilon _i dV} } \right| \leqslant C_1
\end{equation}
and
\begin{equation}\label{eq6.34}
\left| {\int_{\partial T} {\upsilon _i dS} } \right| \leqslant
C_2.
\end{equation}
Since  the  inequality
\begin{equation}\label{eq6.35}
\int_T {\phi ^2 dV \leqslant C} \left( {\int_T {\left| {\nabla
\phi} \right|^2 dV}  +\left| \frac{1} {{Vol\left( T \right)}}
{\int_T {\phi dV} } \right|^2 } \right)
\end{equation}
holds for any $\phi \in H_{1,G}^2 (T)$, taking into account that
(\ref{eq6.27}) and (\ref{eq6.28}) also hold, it follows  that $\{
\upsilon _i \} _{i \in \mathbb{N}}, \upsilon _i \in A$ is bounded
in  $L^2_G (T)$. Moreover, since (\ref{eq6.28}) holds we conclude
that $\mathop {\sup }_{i \in \mathbb{N}} \left( {\left\| {\upsilon
_i } \right\|_{H_{1, G}^2 } } \right) < \infty $. Hence, as in the
previous case there exists $
\upsilon  \in A$ such that $I(\upsilon ) = \mu $.\\
Recall that, if $\nu$ is the Lagrange multiplier, the Euler
equation is
\begin{equation}\label{eq6.37}
\int_T {\nabla ^i \upsilon \nabla _i h} dV + a\int_T h dV +
b\int_{\partial T} {hdS = } \nu \left( {\int_T {fe^\upsilon  } dV
+ \int_{\partial T} {ge^\upsilon  } dS} \right),
\end{equation}
for all $ h \in H_1^2 (T)$. \\
For $h=1$ since $ K(\upsilon ) = 0$ we find
\[
\nu  = -\left( {aVol(T) + bVol(\partial T)} \right)\left( {\int_T
{fe^\upsilon  } dV + \int_{\partial T} {ge^\upsilon  } dS}
\right)^{ - 1}=1.
\]
Using the same arguments as in case $\bf{{1}},\,$ we prove that $
\upsilon  \in C^\infty_G (\overline T)$ and that is a solution of $\;{\bf(P_2)}$.\\

If $g \equiv 0$ we have
\begin{eqnarray*}
  R &=& \left| {\int_T {fe^\upsilon  } dV} \right| \leqslant \left( {\mathop
  {\max }_T \left| f \right|} \right)\int_T {e^\upsilon  } dV \hfill \\
& \leqslant & \left( {\mathop {\max }\limits_T \left| f \right|}
   \right)C_\varepsilon  \exp \left[ {\left( {1 + \varepsilon } \right)\frac{1}
{{16\pi ^2 \left( {l - r} \right)}}\left\| {\nabla \upsilon }
\right\|_2^2  + \frac{1} {{Vol(T)}}\int_T {\upsilon dV} } \right].
\end{eqnarray*}
Hence, if $R < 8\pi ^2 (l - r)$, following the same process as
above
we prove that $\;{\bf(P_2)}$ has a solution.\\

$\bf{3.}\,$ Suppose that $\;R>0$ and $a$, $b$ not both $\geq 0$
(the case $\;R<0$ and $a$, $b$ not both $\leq 0$ can be
treated in the same way).\\
By (\ref{eq6.13}) it is necessary that $\;f$,$\;g$ are not both
$\geq 0$ everywhere. Then $A\not=\emptyset$.

$\bf{(a)}$ $\;a<0$,$\;b>0$,$\;f<0$,$\;g\leq0$ and $\;bVol(\partial
T)<4\pi^2(l-r)$ if $\;g\not\equiv0$ or $\;b
\,Vol(\partial T)<8\pi^2(l-r)$ if $\;g\equiv0$.\\
Since $f \in C^\infty_G  (\overline T)$ is negative everywhere and
$\overline T$ is compact, there exists $\delta  > 0$ such that
$\left| f \right| \geqslant \delta  > 0$.\\
If $\upsilon  \in A$ we have
\begin{equation}\label{eq6.38}
\left| R \right| = \left| {\int_T {fe^\upsilon  } dV +
\int_{\partial T} g e^\upsilon  dS} \right| = \int_T {\left| f
\right|e^\upsilon  } dV + \int_{\partial T} {\left| g \right|}
e^\upsilon  dS
\end{equation}
and by elementary inequality $e^x   \geqslant 1 + x,\;
x\in\mathbb{R} $ we obtain
\[
\left| R \right| \geqslant \int_T {\left| f \right|e^\upsilon  }
dV \geqslant \delta \int_T {e^\upsilon } dV \geqslant \delta
\int_T {(1 + \upsilon )} dV = \delta Vol(T) + \delta \int_T
\upsilon  dV.
\]
Since $a < 0$ we finally obtain
\begin{equation}\label{eq6.39}
a\int_T \upsilon  dV \geqslant a\left( {\frac{{\left| R \right|}}
{\delta } - Vol(T)} \right).
\end{equation}
By (\ref{eq6.23}) implies  that for any $\varepsilon>0$ there
exists a constant $\tilde C_\varepsilon  $ such that
\begin{equation}\label{eq6.40}
\left| R \right| \leqslant \tilde C_\varepsilon  \exp \left[
{\left( {1 + \varepsilon } \right)\frac{1} {{8\pi ^2 \left( {l -
r} \right)}}\left\| {\nabla \upsilon } \right\|_2^2  + \frac{1}
{{Vol(\partial T)}}\int_{\partial T} {\upsilon dS} } \right].
\end{equation}
By (\ref{eq6.40}) we obtain
\begin{equation}\label{eq6.41}
b\int_{\partial T} {\upsilon dS}  \geqslant bVol(\partial T)\ln
{\frac{{\left| R \right|}} {{\tilde C_\varepsilon }}}  - \left( {1
+ \varepsilon } \right)\frac{{bVol(\partial T)}} {{8\pi ^2 \left(
{l - r} \right)}}\left\| {\nabla \upsilon } \right\|_2^2.
\end{equation}
By the definition of $I(\upsilon)$ and (\ref{eq6.39}),
(\ref{eq6.41}) we obtain
\begin{eqnarray}\label{eq6.42}
I(\upsilon )& \geqslant & \left[ {\frac{1} {2} - \left( {1 +
\varepsilon } \right)\frac{{bVol(\partial T)}}
{{8\pi ^2 \left( {l - r} \right)}}} \right]\left\| {\nabla \upsilon } \right\|_2^2 \nonumber \hfill \\
&& + b \;V\!ol(\partial T)\ln {\frac{{\left| R \right|}} {{\tilde
C_\varepsilon  }}}  + a\left( {\frac{{\left| R \right|}} {\delta }
- Vol(T)} \right).
\end{eqnarray}
If $bV\!ol(\partial T) < 4\pi ^2 \left( {l - r} \right)$, that is
$b< \frac{l - r}{lr}$ and $\varepsilon $ is chosen small enough,
(\ref{eq6.42}) implies that $I(\upsilon )$ is bounded bellow for
all $\upsilon  \in A$ and we can prove the existence of a solution
of $\;{\bf(P_2)}$ as in the
previous cases.\\
If $g \equiv 0$, it suffices to assume that $b< \frac{2(l -
r)}{lr}$ and then by (\ref{eq6.40}) we obtain
\[
\left| R \right| = \int_T {\left| f \right|e^\upsilon  } dV
\leqslant \tilde C_\varepsilon  \exp \left[ {\left( {1 +
\varepsilon } \right)\frac{1} {{16\pi ^2 \left( {l - r}
\right)}}\left\| {\nabla \upsilon } \right\|_2^2  + \frac{1}
{{Vol(\partial T)}}\int_{\partial T} {\upsilon dS} } \right]
\]
and we continue as above.

$\bf{(b)}\,$ $\;a>0$,$\;b<0$,$\;f\leq0$,$\;g<0$ and
$\;aVol(T)<4\pi^2(l-r)$.\\
We work as in the previous case and, supposing that $a<
\frac{2(l-r)}{lr^2}$ we conclude the existence of a solution
of $\;{\bf(P_2)}$.\\
Cases {\bf(c)} and {\bf(d)} are similar to {\bf(b)}.\\

$\bf{4.}\,$ Case $\;a\leq0$,$\;b\leq0$, not both $=0$.\\
By  (\ref{eq6.14}) it is necessary to assume that $\int_T f
dV+\int_{\partial T} g dS>0$ and by (\ref{eq6.13}) arises that
$f$, $g$ are not both $ \leqslant 0$ a.e..\\
The proof of this case is based upon the  method of upper
solutions and lower solutions and is the same as the one in
Theorem 2, case $\mathrm{(iv)}$ of \cite{Che4}.\\
Let us sketch the proof: It suffices to find functions $\upsilon _
- ,\upsilon _ + \in C^\infty_G  (\overline T)$ such that $\upsilon
_ + \geqslant \upsilon _ - $ which satisfy the equations
\begin{equation}\label{eq6.43}
\left. \begin{gathered} \Delta \upsilon _ +   + a + fe^{\upsilon _
+  }  \geqslant 0\quad\mathrm{in}\quad
{\rm T}\,\, \hfill \\
\frac{{\partial \upsilon _ +  }} {{\partial n}} + b + ge^{\upsilon
_ +  }  \geqslant 0\quad\mathrm{on}\quad
\partial {\rm T} \hfill \\
\end{gathered}  \right\}
\end{equation}
and
\begin{equation}\label{eq6.44}
\left. \begin{gathered} \Delta \upsilon _ -   + a + fe^{\upsilon _
-  }  \leqslant 0\quad\mathrm{in}\quad
{\rm T}\,\, \hfill \\
\frac{{\partial \upsilon _ -  }} {{\partial n}} + b + ge^{\upsilon
_ -  }  \leqslant 0\quad\mathrm{on}\quad
\partial {\rm T} \hfill \\
\end{gathered}  \right\}
\end{equation}
respectively.\\
We denote by $P_{(a,b)}$ the nonlinear problem $\bf(P_2)$
and solve this case in four steps. More precisely, we prove that: \\
{1}. For any $u \in C ^ 0_G (\overline T)$, $P_{(a,b)} $ accepts a
lower solution $\upsilon _ -  $ such
that $\upsilon _ -   \leqslant u$.\\
{2}. For any $u \in C^0_G (\overline T)$, $P_{(a,b)} $ accepts an
upper solution $\upsilon _ +  $ such that
$\upsilon _ +   \geq u$.\\
{3}. Choosing $f$, $g$ appropriately,
 the set $S_{f,g} $ can be contained in $\mathbb{R}_ - ^2  = \left\{
{\left( {a,b} \right) \ne \left( {0,0} \right):a \leqslant 0,b
\leqslant 0} \right\}$ strictly.\\
{4}. If $f$, $g$ are $\not  \equiv 0$ and nonnegative everywhere
then $S_{f,g}  = \mathbb{R}_ - ^2 $. \mbox{ }\hfill $\Box$

\section{Proofs of the Lemmas}

{{\bf{Proof of Lemma \ref{L2.1}}} $\bf{\;1.}$  Let
$\varepsilon_0>0$ and $(T_j)_{j=1,...,N}$ be a finite covering of
$\bar{T}$, where
$$
T_j = \left\{ {Q \in \mathbb{R}^3 :d(Q,O_{P_j } ) < \delta _j
,\quad\delta _j=l_j\varepsilon_j}\quad\mathrm{and}
\quad\varepsilon _j \leqslant \varepsilon_0\right\}
$$
Then for any $\upsilon \in C_{0,G}^\infty (T_j )$ by (\ref{eq2.1})
we obtain
\begin{eqnarray*}
  \int_{T_j } {e^\upsilon  dV} & =&
  \int_{I \times D} {e^{\upsilon  \circ \xi _j^{ - 1} }
  \left( {\sqrt g  \circ \xi _j^{ - 1} } \right)d\omega dtds}  \hfill \\
  &=& 2\pi l_j \delta _j^2
  \int_D {e^\phi  } (1 + \frac{{\delta _j }}
  {{l_j }}t)dtds \hfill \\
  & \leqslant&   2\pi \,l_j \delta _j^2 (1 + \varepsilon_0 )\int_D {e^\phi  dtds}.  \hfill \\
\end{eqnarray*}
From this and by Theorem $1$ of  \cite{Che1} we have
\[
\int_{T_j } {e^\upsilon  dV \leqslant 2\pi \,l_j \delta _j^2 C(1 +
\varepsilon_0 )\,} \exp \left( {\mu _2 \int\limits_D {\left|
{\nabla \phi } \right|} ^2 dtds} \right),
\]
where $\mu _2=\frac{1}{16\pi}$ is the best constant of Sobolev
inequality
\[
\int_D {e^f dV \leqslant C\exp \left[ {\mu _2 \left\| {\nabla f}
\right\|_2^2 } \right]},
\]
with  $ f \in \stackrel{\circ}{H}\!_1^2 \,(D)$. \\\\
Moreover from (\ref{eq2.2}) we obtain
\begin{eqnarray*}
  \int_{T_j } {\left| {\nabla \upsilon } \right|} ^2 dV &=&
  2\pi \int_D {\left| {\nabla \phi } \right|^2 \left( {l_j  +
  \delta _j t} \right)dtds}  \hfill \\
  &\geqslant &2\pi \,
  l_j \left( {1 - \varepsilon_0 } \right)\int_D {\left| {\nabla \phi } \right|^2 dtds},  \hfill \\
\end{eqnarray*}
thus
\[
\int_D {\left| {\nabla \phi } \right|} ^2 dtds \leqslant \frac{1}
{{2\pi \,l_j }}  \frac{1} {{1 - \varepsilon_0 }}\int_{T_j }
{\left| {\nabla \upsilon } \right|^2 dV}.
\]
Finally, we have
\[
\int_{T_j } {e^\upsilon  dV \leqslant } \,C_{\varepsilon_0} \exp
\left[ {\left( {1 + c\varepsilon_0 } \right)\frac{{\mu _2 }}
{{2\pi \,l_j }}\int_{T_j } {\left| {\nabla \upsilon } \right|} ^2
dtds} \right],
\]
where $ C_{\varepsilon _0 } = 2\pi \,l_l \delta _j^2 C\left( {1 +
\varepsilon_0 } \right)\,$ and $ \frac{1} {{1 - \varepsilon_0 }} =
1 + c\varepsilon_0 \,,c>0$.\\

$\bf{2.}$ Let us choose $\delta>0$ such that the torus $\bar{T}$
is covered by $N$ open subsets
$$
T_{j,\delta/2}=\left\{Q \in T:d\left( Q,O_{P_j }\right) < \delta
/2\right\}
$$
We consider  the decreasing real valued $C^\infty$ function
$\Psi(r)$ , which equals $1$ for $0\leq r\leq \delta/2$ and $0$
for $r\geq\delta$ and we note $ \Psi _j (Q) = \Psi(d(Q,O_{P_j }
))$.\\
The $\Psi_j$'s defined on $T_j=\left\{Q \in T:d\left( Q,O_{P_j
}\right) < \delta \right\}$ are $G-$invariant, but they are not a
partition of
unity.\\
Let  $ \upsilon  \in C_{0,G}^\infty  (T)$. Then $ (\upsilon \Psi
_j )\in C_{0,G}^\infty  (T_j )$ and from the first part of this
lemma we obtain
\begin{equation}\label{eq4.25}
\int_{T_j } {e^{\upsilon \Psi _j } dV \leqslant } \,C\exp \left[
{\left( {1 + c\varepsilon_0 } \right)\frac{{\mu _2 }} {{2\pi l_j
}}\left\| {\nabla \left( {\upsilon \Psi _j } \right)} \right\|_2^2
} \right]
\end{equation}
Because of the following relation
\[
\begin{gathered}
  \left\| {\nabla \left( {\upsilon \Psi _j } \right)} \right\|_2^2  \leqslant
  \left\| {\Psi _j \nabla \upsilon } \right\|_2^2  +
  2\left\| {\Psi _j \nabla \upsilon } \right\|\left\| {\upsilon \nabla \Psi _j } \right\| +
  \left\| {\upsilon \nabla \Psi _j } \right\|_2^2 . \hfill \\
\end{gathered}
\]
and since for all $ \varepsilon _0  > 0$ a constant
$D_{\varepsilon _0 }$ exists such that
\[ \left\| {\Psi _j \nabla
\upsilon } \right\|\left\| {\upsilon \nabla \Psi _j } \right\|
\leqslant \varepsilon _0 \left\| {\Psi _j \nabla \upsilon }
\right\|_2^2  + D_{\varepsilon _0 } \left\| {\upsilon \nabla \Psi
_j } \right\|_2^2,
\]
we obtain
\[
\begin{gathered}  \left\| {\nabla \left( {\upsilon \Psi _j }
 \right)} \right\|_2^2  \leqslant
  \left( {1 + 2\varepsilon _0 } \right)\left\| {\nabla \upsilon } \right\|_2^2  +
  \tilde D\left\| \upsilon  \right\|_2^2,  \hfill \\
\end{gathered}
\]
where $ \tilde D = \left( {2D_{\varepsilon _0 }  + 1}
\right)\left( {\sup _T \,\left| {\nabla \Psi _j } \right|^2 }
\right)$.\\
From (\ref{eq4.25}) because of the last inequality we have
\begin{equation}\label{eq4.26}
\int_{T_j } {e^{\upsilon \Psi _j } dV \leqslant } \,C\exp \left[
{\left( {1 + c\varepsilon_0 } \right)\left( {1 + 2\varepsilon _0
}\right)\frac{{\mu _2 }} {{2\pi \,l_j }} \left\| {\nabla \upsilon
} \right\|_{^2 }^2 + D\left\| \upsilon \right\|_2^2 } \right],
\end{equation}
where $ D = \left( {1 + c\varepsilon_0 } \right)\frac{{\mu _2 }}
{{2\pi \,l_j }}\tilde D$.\\
 Since $\inf l_j  = l - r$ given $\varepsilon>0$ we can choose
 $\varepsilon_0$ small enough such that from (\ref{eq4.26})  we obtain
\begin{eqnarray*}
\int_T {e^\upsilon  } dV &\leqslant &\sum\limits_{i = 1}^N {\int_{T_{j,\delta /2} } {e^\upsilon  } dV}\hfill \\
  & \leqslant& \sum\limits_{i = 1}^N {\int_{T_{j,\delta } } {e^{\upsilon \Psi _j } } dV}  \hfill \\
  &\leqslant& C\exp \left[ {\left( {\frac{{\mu _2 }}
{{2\pi \left( {l - r} \right)}} + \varepsilon } \right)\left\|
{\nabla \upsilon } \right\|_2^2  + D\left\| \upsilon \right\|_2^2
} \right],
\end{eqnarray*}
and so  we have the desired inequality.\\
Now we need to prove that the constant $ \frac{{\mu _2 }} {{2\pi
(l - r)}}$  is the best constant $\mu$ such that the inequality
\[
\int_T {e^\upsilon  dV \leqslant C\exp \left[ {\left( {\mu +
\varepsilon } \right)\left\| {\nabla \upsilon } \right\|_2^2  +
D\left\| \upsilon  \right\|_2^2 } \right]}
\]
holds for all $\upsilon\in \stackrel{\circ}{H}\!_{1,G}^2$ .\\
For that purpose, for all $\varepsilon$, we need to find a
sequence $(\upsilon_\alpha)\in \stackrel{\circ}{H}\!_{1,G}^2$,
such that for all $ \Delta ,{\rm E} \in \mathbb{R}$ the following
holds:
\begin {equation}\label{eq4.27}
\mathop {\lim }\limits_{\alpha  \to 0} \frac{{\left\| {\nabla
\upsilon _\alpha  } \right\|_2^2  + \Delta \left\| {\upsilon
_\alpha  } \right\|_2^2  + {\rm E}}} {{\ln \int_T {e^{\upsilon
_\alpha  } dV} }} \leqslant \frac{{2\pi (l - r)}} {{\mu _2 }} +
\varepsilon.
\end {equation}
Let us consider the orbit $O_{inf}$ of minimum length $2\pi(l-r)$.
For any $\varepsilon_0>0$, let
$$
T_{j_{_0}} = \left\{ {Q \in \mathbb{R}^3:d\left(Q,O_{inf} \right)
< \delta },\quad \delta=\varepsilon_0 (l-r)\right\} ,
$$
where $d(Q,O_{inf})$ denotes the distance from $Q$ to the orbit
$O_{inf}$.\\
It is easy to prove that
\begin{eqnarray}\label{eq4.120}
 d(Q,O_P ) = \delta d_D (\xi_{j_{_0}}
(Q),O ) = \delta \sqrt {t^2  + s^2 },
\end{eqnarray}
where $d_D$ denotes the distance in the disc $D$ centered on
$O$.\\
For all $\alpha>0$ define the functions $(\upsilon_\alpha)$ by
$$
\upsilon_\alpha(Q)=\left\{
\begin{array}{c}
  -2ln\left(\alpha+d^2(Q,O_{inf})\right)+2ln\left(\alpha+\delta^2\right),
  \mathrm{if}\;Q\in T \cap T_{j_{_0}}\\
  0\;\;\;\;\;\;\;\;\;\;\;\;\;\;\;\;\;\;\;\;\;\;\;\;\;\;\;\;\;\;\;\;\;\;\;\;\;\;\;\;\;\;\;\;\;\;\;\;\;\;\;\;\;\;\;
\;\;\;\;\;  ,\mathrm{if}\;Q \in T \,\backslash \,T_{j_{_0}} \\
\end{array}\right.
$$
Since $\upsilon_\alpha$ depends only on the distance to $O_{inf}$,
$\upsilon_\alpha\in\stackrel{\circ}{H}\!_{1,G}^2(T_{j_{_0}})$.\\
Setting $\phi_\alpha=\upsilon_\alpha\circ\xi_{j_{_0}}^{-1}$  we
obtain
\begin{eqnarray}\label{eq4.28}
  \int_T {e^{\upsilon _\alpha  } dV}  &=&
   \int_{I \times D} {e^{\upsilon _\alpha   \circ \xi _{j_{_0}} ^{ - 1} } }
   (\sqrt g  \circ \xi_{j_{_0}} ^{ - 1} )d\omega dtds \nonumber\hfill \\
    & =& \int_{I \times D}
  {e^{\phi _\alpha  } }
   \delta ^2 \bigl((l-r)  + \delta t\bigl )d\omega dtds \nonumber\hfill \\
  &=& 2\pi \,(l-r) \delta ^2
   \int_D {e^{\phi _\alpha  } }
  \left (1 + \frac{\delta }
{{l-r}}t\right )dtds \nonumber\hfill \\
   & \geq &2\pi \,(l-r) \delta ^2
  (1 - \varepsilon_0 )\int_D {e^{\phi _\alpha  } } dtds.
\end{eqnarray}
Hence, by definition of $\upsilon_\alpha$ and because of
(\ref{eq4.120})  for all $ \xi _{j_{_0}} (Q) = (t,s) \in D$ we
obtain
\[
\phi _\alpha  (\xi _j (Q)) =  \ln \left( {\frac{{\alpha  + \delta
^2 }} {{\alpha  + \delta ^2 \left( {t^2  + s^2 } \right)}}}
\right)^2,
\]
thus
\[ \int_D {e^{\phi _\alpha  } } dtds = \int_D
{\left( {\frac{{\alpha  + \delta ^2 }} {{\alpha  + \delta ^2
\left( {t^2  + s^2 } \right)}}} \right)^2 } dtds.
\]
Changing variables in the latter equality we obtain
\begin{eqnarray}\label{eq4.29}
  \int_D {e^{\phi _\alpha  } } dtds & = & \int_0^{2\pi } {\int_0^1
  {\frac{{\left( {\alpha  + \delta ^2 } \right)^2 r}}
{{\left( {\alpha  + \delta ^2 r^2 } \right)^2 }}} } drd\theta
\,\nonumber \\ &=& \frac{{\left( {\alpha  + \delta ^2 } \right)^2
\pi }} {{\delta ^2 }}\int_0^1 {\frac{{\left( {\alpha  + \delta ^2
r^2 } \right)^\prime  }} {{\left( {\alpha  + \delta ^2 r^2 }
\right)^2 }}} dr \nonumber \\ &=& \frac{{\left( {\alpha  + \delta
^2 } \right)\pi }} {\alpha }.
\end{eqnarray}
By (\ref{eq4.28}) and (\ref{eq4.29}) we have
\begin{eqnarray*}
\int_T {e^{\upsilon _\alpha  } dV} \, &\geq&
  \left( {1 - \varepsilon_0 } \right)2\pi \,(l-r) \delta ^2
  \frac{{\left( {\alpha  + \delta ^2 } \right)\pi }}
{\alpha } \hfill \\
 & \geqslant&
   \left( {1 - \varepsilon_0 } \right)2\pi ^2 \,(l-r) \delta ^4 \frac{1}
{\alpha } \hfill \\
& =& C_{\varepsilon_0}  \frac{1} {\alpha },
\end{eqnarray*}
and then
\begin{equation}\label{eq4.30}
\ln \int_T {e^{\upsilon _\alpha  } dV} \, \geqslant \ln
C_{\varepsilon_0}   + \ln \frac{1} {\alpha },
\end{equation}
where $C_{\varepsilon_0} = \left( {1 - \varepsilon_0 } \right)2\pi
^2
\,(l-r) \delta ^4$.\\
Moreover, because of (\ref{eq2.2}) we have
\begin{eqnarray}\label{eq4.31}
  \left\| {\nabla \upsilon _\alpha  } \right\|_2^2  &=&
  2\pi \int_D {\left| {\nabla \phi _\alpha
  \left( {t,s} \right)} \right|} ^2 \bigl( {\left( {\,l - r} \right) + \delta t} \bigl)dtds \nonumber \\
   & \leqslant& \left( {1 + \varepsilon_0 } \right)2\pi \,\left( {\,l - r}
    \right)\int_D {\left| {\nabla \phi _\alpha  \left( {t,s} \right)} \right|} ^2
    dtds.
\end{eqnarray}
Since
\begin{eqnarray*}
  \left| {\nabla \phi _\alpha  \left( {t,s} \right)} \right|^2 &
  =&
  \left| {\nabla \left[ { - 2\ln \left( {\alpha  + \delta ^2 \left( {t^2
   + s^2 } \right)} \right) + 2\ln \left( {\,\alpha  + \delta ^2 } \right)}
    \right]} \right|^2  \\
&=&\left| { - 2\nabla \left[ {\ln \left( {\alpha  + \delta ^2
\left( {t^2
      + s^2 } \right)} \right)} \right]} \right|^2 \\
 &=&
  4\left| {\left( {\frac{{2\delta ^2 t}}
{{\alpha  + \delta ^2 \left( {t^2  + s^2 }
\right)}},\frac{{2\delta ^2 s}} {{\alpha  + \delta ^2 \left( {t^2
+ s^2 } \right)}}} \right)} \right|^2   \\
 & = &\frac{{16\delta ^4 \left( {t^2  + s^2 } \right)}}
{{\left[ {\alpha  + \delta ^2 \left( {t^2  + s^2 } \right)}
\right]^2 }},
\end{eqnarray*}
we have
\begin{eqnarray*}
  \int_D {\left| {\nabla \phi _\alpha  (t,s)} \right|} ^2 dtds
  &=&   \int_D {\frac{{16\delta ^4 \left( {t^2  + s^2 } \right)}}
{{\left[ {\alpha  + \delta ^2 \left( {t^2  + s^2 } \right)} \right]^2 }}dtds} \\
&=&  2\pi \int_0^1 {\frac{{16\delta ^4 r^2 }} {{\left( {\alpha  +
\delta ^2 r^2 } \right)^2 }}}r dr.
\end{eqnarray*}
Changing variables we obtain
\begin{eqnarray}\label{eq4.32}
\int_D {\left| {\nabla \phi _\alpha  (t,s)} \right|} ^2 dtds &=&
16\pi \int_0^{\delta ^2 } {\frac{\tau } {{\left( {\alpha  + \tau }
\right)^2 }}} d\tau \nonumber\\ & =& \frac{1} {{\mu _2
}}\int_0^{\delta ^2 } {\frac{\tau } {{\left( {\alpha  + \tau }
\right)^2 }}} d\tau.
\end{eqnarray}
We further define the function $$ h(\alpha ) = \ln \frac{1}
{\alpha } - \int_0^{\delta ^2 } {\frac{\tau } {{\left( {\alpha  +
\tau } \right)^2 }}} d\tau ,\alpha  > 0$$
 and changing the variable we obtain
\begin{eqnarray*}
h(\alpha )& =&  \ln \frac{1} {\alpha } - \int_0^{\delta ^2 }
{\frac{\tau } {{\alpha ^2 \left( {1 + \frac{\tau } {\alpha }}
\right)^2 }}} d\tau \\& =& \int_{\delta ^2 }^{\delta ^2 /\alpha }
{\frac{1} {u}} du - \int_0^{\delta ^2 /\alpha } {\frac{u}
{{\left( {1 + u} \right)^2 }}} du  \\
&=& \int_{\delta ^2 }^{\delta ^2 /\alpha } {\left( {\frac{1} {u} -
\frac{u} {{\left( {1 + u} \right)^2 }}} \right)} du -
\int_0^{\delta ^2 } {\frac{u}
{{\left( {1 + u} \right)^2 }}} du \\
& =& \int_{\delta ^2 }^{\delta ^2 /\alpha } {\left[ {1 - \left( {1
+ u^{ - 1} } \right)^{ - 2} } \right]u^{ - 1} } du -
\int_0^{\delta ^2 } {\frac{u} {{\left( {1 + u} \right)^2 }}} du
\end{eqnarray*}
and because of
\begin{eqnarray*}
\left[ {1 - \left( {1 + u^{ - 1} } \right)^{ - 2} } \right]u^{ -
1} & =& \left[ {1 - \left( {1 + \frac{1} {u}} \right)^{ - 2} }
\right]\frac{1} {u} = \left[ {1 - \left( {\frac{u} {{u + 1}}}
\right)^2 } \right]\frac{1}
{u}  \\
& =& \frac{1} {{u + 1}} \cdot \frac{{2u + 1}} {{u + 1}} \cdot
\frac{1} {u} < \frac{2} {{u + 1}} \cdot \frac{1} {u} < \frac{2}
{{u^2 }}
\end{eqnarray*}
we finally obtain
\begin{eqnarray*}
\mathop {\lim }\limits_{\alpha  \to 0} h(\alpha )& =& \int_{\delta
^2 }^\infty
  {\left[ {1 - \left( {1 + u^{ - 1} } \right)^{ - 2} } \right]u^{ - 1} } du - \int_0^{\delta ^2 }
  {\frac{u}
{{\left( {1 + u} \right)^2 }}} du  \\
&\leq &\int_{\delta ^2 }^\infty  {2u^{ - 2} } du -
  \int_0^{\delta ^2 } {\frac{u}
{{\left( {1 + u} \right)^2 }}} \,du = C_0.
\end{eqnarray*}
Thus, for any $\alpha>0$ close to $0$ the following holds:
\begin{equation}\label{eq4.33}
\int_0^{\delta ^2 } {\frac{\tau } {{\left( {\alpha  + \tau }
\right)^2 }}} d\tau  = \ln \frac{1} {\alpha } + C_1.
\end{equation}
From (\ref{eq4.31}), (\ref{eq4.32}) and (\ref{eq4.33}) we obtain
\begin{equation}\label{eq4.34}
\left\| {\nabla \upsilon _\alpha  } \right\|_2^2  \leqslant
\frac{{\left( {1 + \varepsilon_0 } \right)2\pi \,(l-r) }} {{\mu _2
}}\ln \frac{1} {\alpha } + C.
\end{equation}
On the other hand we have
\begin{eqnarray*}
\left\| {\upsilon _\alpha  } \right\|_2^2 & =& 2\pi \delta ^2
\int_D {\left| {\phi _\alpha  } \right|} ^2 \left( {\left( {l - r}
 \right) + \delta t} \right)dtds\, \\
& \leqslant &(1 + \varepsilon _0 )2\pi \left(
 {l - r} \right)\delta ^2 \int_D {\left| {\phi _\alpha  }
  \right|} ^2 dtds \\
& =& C_0 \int_D {\left| {2\ln \left(
  {\frac{{\alpha  + \delta ^2 }}
{{\,\alpha  + \delta ^2 \left( {t^2  + s^2 } \right)}}} \right)}
\right|}
 ^2 dtds  \\
& =& 4C_0 \int_D {\left( {\ln
   \left( {\alpha  + \delta ^2 } \right) - \ln \left( {\alpha  + \delta ^2
   \left( {t^2  + s^2 } \right)} \right)} \right)} ^2 dtds\, \\
&\leqslant& 8\pi C_0 \int_D {\left[ {2\ln ^2
   \left( {\alpha  + \delta ^2 } \right) + 2\ln ^2 \left( {\alpha  +
    \delta ^2 \left( {t^2  + s^2 } \right)} \right)} \right]} dtds\,  \\
& =& 8\pi C_0 \left( {\int_0^1 {2\ln ^2 \left( {\alpha  + \delta
^2 } \right)} rdr +
 \int_0^1 {2\ln ^2 \left( {\alpha  + \delta ^2 r^2 } \right)} \,rdr} \right) \hfill \\
& =& 8\pi C_0 \ln ^2 \left( {\alpha  + \delta ^2 } \right)\int_0^1
2 rdr + \frac{{8\pi C_0 }}
{{\delta ^2 }}\int_0^1 {\ln ^2 \left( {\alpha  + \delta ^2 r^2 } \right)} 2\delta ^2 rdr \hfill \\
& = &8\pi C_0 \ln ^2 \left( {\alpha  + \delta ^2 } \right) +
\frac{{8\pi C_0 }} {{\delta ^2 }}\int_0^1 {\ln ^2 \left( {\alpha
+ \delta ^2 r^2 } \right)}
\left( {\alpha  + \delta ^2 r^2 } \right)^\prime  dr \hfill \\
& =& C_1  + C_2 \int_\alpha ^{\alpha  + \delta ^2 } {\ln ^2 \zeta } d\zeta  \hfill \\
\end{eqnarray*}
\begin{eqnarray*}
& =& C_1  + C_2 \left[ {\zeta \left( {\ln ^2 \zeta  - 2\ln \zeta
+ 2}
\right)} \right]_\alpha ^{\alpha  + \delta ^2 }  \hfill \\
& =& C_1  + C_2 \left[ {\left( {\alpha  + \delta ^2 } \right)
\left( {\ln ^2 \left( {\alpha  + \delta ^2 } \right) - 2\ln \left(
{\alpha  + \delta ^2 } \right) + 2}
\right)} \right] \hfill \\
&& - C_2 \alpha \left( {\ln ^2 \alpha  - 2\ln \alpha  + 2}
\right),
\end{eqnarray*}
and since $\mathop {\lim }\limits_{\alpha  \to 0^{^+ }} \left(
{\alpha \ln \alpha } \right) = \mathop {\lim }\limits_{\alpha \to
0^{^+ }} \left( {\alpha \ln ^2 \alpha } \right) = 0$ we have
\begin{equation}\label{eq4.35}
\left\| {\upsilon _\alpha  } \right\|_2^2  \leqslant \,\,C_1  +
C_2 C_3  = C.
\end{equation}
Finally, from (\ref{eq4.30}), (\ref{eq4.34}) and (\ref{eq4.35})
for any $ \Delta ,{\rm E} \in \mathbb{R}$ the following holds:
\[
\frac{{\left\| {\nabla \upsilon _\alpha  } \right\|_2^2  + \Delta
\left\| {\upsilon _\alpha  } \right\|_2^2  + {\rm E}}} {{\ln
\int_T {e^{\upsilon _\alpha  } dV} }} \leqslant
\frac{{\frac{{\left( {1 + \varepsilon _0 } \right)2\pi \left( {l -
r} \right)}} {{\mu _2 }}\ln \frac{1} {\alpha } + C}} {{\ln
\frac{1} {\alpha } + \ln C_{\varepsilon_{_0}}}},
\]
thus,
\begin{equation}\label{eq4.36}
\mathop {\lim }\limits_{\alpha  \to 0} \frac{{\left\| {\nabla
\upsilon _\alpha  } \right\|_2^2  + \Delta \left\| {\upsilon
_\alpha  } \right\|_2^2  + {\rm E}}} {{\ln \int_T {e^{\upsilon
_\alpha  } dV} }} \leqslant \left( {1 + \varepsilon _0 }
\right)\frac{{2\pi \left( {l - r} \right)}} {{\mu _2 }}.
\end{equation}
For any $ \varepsilon  > 0$ consider $ \varepsilon_0  > 0$ such
that $ \left( {1 + \varepsilon _0 } \right)\frac{{2\pi \left( {l -
r} \right)}} {{\mu _2 }} \leqslant \frac{{2\pi \left( {l - r}
\right)}} {{\mu _2 }} + \varepsilon$ and so from (\ref{eq4.36}) we
obtain our result.\mbox{ }\hfill $\Box$\\
\\
{{\bf{Proof of Lemma \ref{L2.2}}} Following arguments similar to
those in \cite{Aub4} and \cite{Fag1} we prove the first and second
part of the lemma, respectively.

$\bf{1.}$$\,\,\,$ Our aim here is to find a constant
$C_\varepsilon$, such that for any $\varepsilon>0$ and for all
functions $\upsilon\in \mathcal {H}_G$, with $\int_ T \upsilon dV
=0$ the following inequality holds
\[
\int_T {e^\upsilon  dV}  \leqslant C_\varepsilon\exp \left[
{\left( {\mu + \varepsilon } \right)\left\| {\nabla \upsilon }
\right\|_2^2 } \right] ,
\]
where $\mu=\frac{1}{16\pi L}$ if
$\mathcal{H}_G=\stackrel{\circ}{H}\!_{1,G}^2$ and
$\mu=\frac{1}{8\pi L}$ if $\mathcal{H}_G=H_{1,G}^2$.\\

$\bf{(i)}$$\,\,\,$ Let $\upsilon  \in C_{0,G}^\infty (T)$ with
$\int_T \upsilon dV =0$ and $\check{\upsilon}=sup (\upsilon,0)$.\\
Then $\check{\upsilon } \in \stackrel{\circ}{H}\!_{1,G}^2 (T) $
and
$$
 \int_{\mathbb{R}^3 } {\check \upsilon dx = } \frac{1}
{2}\int_{\mathbb{R}^3 } {\upsilon dx},\quad  \int_{\mathbb{R}^3 }
{\left| {\nabla \check \upsilon } \right|dx \leqslant }
\int_{\mathbb{R}^3 } {\left| {\nabla \upsilon } \right|dx}.
$$
For any $ t \in \mathbb{R} $, denote by  $ m_t (\upsilon )$ the
measure of the set
$$
\Omega _t(\upsilon )=\{ x \in T:\upsilon (x) \geqslant t\}.
$$
Given $ \upsilon  \in C_{0,G}^\infty (T)$,  $m_t (\upsilon )$ is a
decreasing function of $t$, not necessarily continuous. Let $m>0$
depending on $\varepsilon$. Then, for a given $ \upsilon \in
C_{0,G}^\infty (T)$ two different cases can occur: whether there
exists $ s \geqslant 0$ such that $ m_s (\upsilon ) \geqslant m $
or not.

$(a)\;$ Suppose there exists $ s \geqslant 0$ such that $
m_s (\upsilon ) \geqslant m $.\\
If we denote
$$
S = \sup \{ s \in \mathbb{R}:m_s (\upsilon ) \geqslant m\},
$$
we will have $ S \geqslant 0$, $ m_{S + 1} (\upsilon ) < m$ and $
m_{S/2} (\upsilon
) \geqslant m$. \\
According to Lemma \ref {L2.1} we have the following
\begin{eqnarray}\label{eq5.13}\
\int_T {e^\upsilon  dV }&=& e^{S + 1}  \int_T {e^{\upsilon - (S +
1)} dV \leqslant } e^{S + 1} \int_T
e^{\check{\widehat{\upsilon  - (S + 1)}}} dV  \nonumber \\
& \leqslant& e^{S + 1} C_{\varepsilon /2} \exp \left[ {\left( {\mu
+ \frac{\varepsilon } {2}} \right)\left\| {\nabla \upsilon }
\right\|_2^2  + D_{\varepsilon /2} \left\|
{\check{\widehat{\upsilon  - (S + 1)}}} \right\|_2^2 } \right]. \nonumber \\
\end{eqnarray}
Since $\int_ T \upsilon dV=0$, and $\left\|
\check{\upsilon}\right\| _1=\frac{1}{2}\left\| \upsilon\right\|
_1$ by Poincar$\grave{e}$ inequality there exists a constant $C_1$
such that
\begin{equation}\label{eq5.14}
\left\| {\check \upsilon } \right\|_1 = \frac{1} {2}\left\|
\upsilon  \right\|_1  \leqslant C_1 \left\| {\nabla \upsilon }
\right\|_2,
\end{equation}
and since $ S + 1 > 0$ we obtain
\begin{equation}\label{eq5.15}
\left\| {\check{\widehat{\upsilon  - (S + 1)}}} \right\|_1
\leqslant \left\| {\bar \upsilon } \right\|_1 \leqslant C_1
\left\| {\nabla \upsilon } \right\|_2.
\end{equation}
From the last two inequalities, by H\"older's inequality and the
Sobolev continuous and compact embedding of
$\stackrel{\circ}{H}\!_{1,G}^2$ in $ L_G^p (T)$, we obtain
\begin{eqnarray}\label{eq5.16}
  \left\| {\check{\widehat{\upsilon  - (S + 1)}}} \right\|_2^2
 & \leqslant& m_{_{S+1} }^{1/2} (\upsilon )\left\|
  {\check{\widehat{\upsilon  - (S + 1)}}} \right\|_4^2  \nonumber \\
& <&  m^{1/2} \left\| {\check{\widehat{\upsilon  - (S + 1)}}} \right\|_4^2   \nonumber \\
&\leqslant & m^{1/2} C_2 \left\| {\nabla \upsilon } \right\|_2^2,
\end{eqnarray}
where $C_2$ is a constant independent of $\upsilon$
and $\mu$.\\
By the definition of $ \Omega _t (\upsilon )$ we have that
$$
 \Omega _0(\upsilon ) = \{ x \in T:\upsilon (x) \geqslant 0\}$$ and $$
\Omega _{S/2} (\upsilon ) = \{ x \in T:\upsilon (x) \geqslant
S/2\}.
$$
Thus
\begin{equation}\label{eq5.17}
\left\| {\check \upsilon } \right\|_1  = \int_{\Omega _0 (\upsilon
)} {\upsilon dV \geqslant } \int_{\Omega _{S/2} (\upsilon )}
{\upsilon dV \geqslant } \int_{\Omega _{S/2} (\upsilon )}
{\frac{S} {2}dV = } \frac{S} {2}m_{S/2} (\upsilon ).
\end{equation}
From (\ref{eq5.14}) and (\ref{eq5.17}), since
$m_{S/2}(\upsilon)\geq m$, we obtain
\begin{equation}\label{eq5.18}
S \leqslant \frac{2} {{m_{S/2} (\upsilon )}}\left\| {\check
\upsilon } \right\|_1  \leqslant \frac{2} {m}C_1 \left\| {\nabla
\upsilon } \right\|_2.
\end{equation}
The elementary inequality
$$
 x < Sx^2  + \frac{1} {S},x \in
\mathbb{R},S > 0,
$$
 with $x = \left\| {\nabla \upsilon }
\right\|_2$ yields
\begin{equation}\label{eq5.19}
\left\| {\nabla \upsilon } \right\|_2  < S\left\| {\nabla \upsilon
} \right\|_2^2  + \frac{1} {S},\,\,S > 0.
\end{equation}
From (\ref{eq5.18}), and because of (\ref{eq5.19}), we obtain
\[
 S\leqslant \frac{2C_1} {m} \left( {S\left\| {\nabla \upsilon }
\right\|_2^2 + \frac{1} {S}} \right),
\]
and with $ \frac {2C_1 }{m }=\frac{m}{S}$ we obtain
\begin{equation}\label{eq5.20}
S \leqslant m\left\| {\nabla \upsilon } \right\|_2^2  + 4C_1^2 m^{
- 3}  = m\left\| {\nabla \upsilon } \right\|_2^2  + C_3 m^{ - 3}.
\end{equation}
Thus, from (\ref{eq5.13}), and because of (\ref{eq5.16}) and
(\ref{eq5.20}), we obtain
\begin{equation}\label{eq5.21}
\int_T {e^\upsilon  dV \leqslant } \,C_\varepsilon  \exp \left[
{\left( {\frac{{\mu _2 }} {{L}} + \frac{\varepsilon } {2} +
D_{\varepsilon /2} C_2 m^{1/2}  + m} \right)\left\| {\nabla
\upsilon } \right\|_2^2 } \right]
\end{equation}
where $ C_\varepsilon   = C_{\varepsilon /2} \exp \left( {C_3 m^{
- 3}  + 1} \right)$.

$(b)\;$ Suppose now that $m_s(\upsilon)<m$ for any $s\geq0$.\\
By Lemma \ref{L2.1} we have the following
\[
\int_T {e^\upsilon  dV \leqslant } \int_T {e^{\check \upsilon }
dV}  \leqslant C_{\varepsilon/2}  exp\left[ {\left( \mu +
\frac{\varepsilon } {2} \right)\left\| {\nabla \check \upsilon }
\right\|_2^2  + D_{\varepsilon/2} \left\| {\,\check \upsilon \,}
\right\|_2^2 } \right]
\]
or
\begin{equation}\label{eq5.22}
\int_T {e^\upsilon  dV \leqslant } C_{\varepsilon/2} \exp \left[
{\left( \mu+ \frac{\varepsilon } {2} \right)\left\| {\nabla
\upsilon } \right\|_2^2  + D_{\varepsilon/2} \left\| {\,\upsilon
\,} \right\|_2^2 } \right].
\end{equation}
In this case, $m_0 (\upsilon)\leq m $ and so
\begin{equation}\label{eq5.23}
 \left\| {\check \upsilon }
\right\|_2^2  \leqslant m_{_0 }^{1/2} (\upsilon )\left\| {\check
\upsilon } \right\|_4^2  < m^{1/2} \left\| {\check \upsilon }
\right\|_4^2  \leqslant m^{1/2} C_2 \left\| {\nabla \upsilon }
\right\|_2^2.
\end{equation}
From (\ref{eq5.22}), and because of (\ref{eq5.23}) we obtain
\begin{equation}\label{eq5.24}
\int_T {e^\upsilon  dV \leqslant } \,C_{\varepsilon/2} \exp \left[
{\left( \mu + \frac{\varepsilon } {2} + D_{\varepsilon/2} C_2
m^{1/2} + m \right)\left\| {\nabla \upsilon } \right\|_2^2 }
\right].
\end{equation}
In both cases we have to choose $m>0$ such that
$$
D_{\varepsilon/2}  C_2 m^{1/2}  + m < \frac{\varepsilon } {2}
$$
and
$$
C_\varepsilon = C_{\varepsilon/2} \exp \left( {\frac{{2C_1
}} {{m^2 }} + 1} \right)
$$
so,  for all $\upsilon\in\stackrel{\circ}{H}\!_{1,G}^2$ with
$\int_T \upsilon dV=0$ the following inequality holds
\begin{equation}\label{eq5.25}
\int_T {e^\upsilon  dV \leqslant C_\varepsilon\exp \left[ {\left(
\mu + \varepsilon  \right)\left\| {\nabla \upsilon } \right\|_2^2
} \right]},
\end{equation}
where $\mu=\frac{1}{16\pi L}$.\\

$\bf{(ii)}$$\,\,\,$  Let now $\upsilon \in H^2_{1,G}$. Following
the same steps as in the  first part of Lemma \ref{L2.1} by
Theorem $3$ of \cite{Che3}, for all $ \upsilon  \in C_G^\infty
(T)$, we obtain
\begin{equation}\label{eq5.26}
\int_{T_j } {e^\upsilon  dV \leqslant } \,C \exp \left[ {(1 +
c\varepsilon )\frac{1 } {{16\pi^2 \,l_j }}\int_{T_j } {\left|
{\nabla \upsilon } \right|} ^2 dV} \right].
\end{equation}
Consequently, since $ C_G^\infty(T)$ is dense in $H^2_{1,G}$ and
(\ref{eq5.26}) holds for any $j=1,2,..,N$, by the second part of
Lemma \ref{L2.1}, we conclude that, for all $\varepsilon>0$, there
are constants $C_\varepsilon$ and $D_\varepsilon$ such that for
all $\upsilon \in H^2_{1,G}$ the following holds
\[
\int_T {e^\upsilon  dV}  \leqslant C_\varepsilon \exp \left[
{\left( \mu + \varepsilon  \right)\left\| {\nabla \upsilon }
\right\|_2^2  + D_\varepsilon\left\| \upsilon  \right\|_2^2 }
\right],
\]
where $\mu=\frac{1}{8\pi L}$ is the best constant for this
inequality.\\
Following the same steps as in the first part of this lemma we
derive that for any $\varepsilon>0$ and for all functions
$\upsilon\in H_{1,G}^2$, with $\int_ T \upsilon dV =0$ the
following inequality holds
\begin{equation}\label{eq5.27}
\int_T {e^\upsilon  dV}  \leqslant C_\varepsilon\exp \left[
{\left( {\mu + \varepsilon } \right)\left\| {\nabla \upsilon }
\right\|_2^2 } \right],
\end{equation}
where $\mu=\frac{1}{8\pi L}$ .\\

By parts $\bf{(i)}$ and $\bf{(ii)}$ of the lemma we conclude that,
for all $\varepsilon>0$, there exists constant $C_\varepsilon$
such that for all $\upsilon\in\stackrel{\circ}{H}\!_{1,G}^2$ or
$\upsilon\in H^2_{1,G}$ with $\int_T \upsilon dV=0$,
inequalities  (\ref{eq5.25}) and (\ref{eq5.27}) hold respectively.\\
Finally, we observe that if $ \tilde \upsilon  = \upsilon -
\frac{1} {2\pi^2r^2l}\int_T {\upsilon dV} $ we have
\begin{eqnarray*}
\int_T \tilde \upsilon  dV&=&\int_T \left(\upsilon - \frac{1}
{2\pi^2r^2l}\int_T {\upsilon dV} \right)dV=\int_T \upsilon dV -
\frac{1} {2\pi^2r^2l}\int_T {\upsilon dV} \int_T dV\\&=&\int_T
\upsilon dV - \frac{1} {Vol(T)}Vol(T)\int_T {\upsilon dV}=0,
\end{eqnarray*}
and so, rewriting (\ref{eq5.25}) and (\ref{eq5.27}) with $ \tilde
\upsilon  = \upsilon - \frac{1} {2\pi^2r^2l}\int_T {\upsilon dV} $
we obtain:
\begin{eqnarray*}
\int_T {e^{\upsilon - \frac{1} {2\pi^2r^2l}\int_T \upsilon dV} dV
\leqslant C_\varepsilon\exp \left[ {\left( \mu + \varepsilon
\right)\left\| {\nabla \left(\upsilon - \frac{1}
{2\pi^2r^2l}\int_T {\upsilon dV}\right) } \right\|_2^2 } \right]},
\end{eqnarray*}
\begin{eqnarray*}
e^{ - \frac{1} {2\pi^2r^2l}\int_T \upsilon dV} \int_T {e^{\upsilon
} dV \leqslant C_\varepsilon\exp \left[ {\left( \mu + \varepsilon
\right)\left\| {\nabla \upsilon  } \right\|_2^2 } \right]}
\end{eqnarray*}
or
\begin{eqnarray*}
 \int_T {e^{\upsilon
} dV \leqslant C_\varepsilon\exp \left[ {\left( \mu + \varepsilon
\right)\left\| {\nabla \upsilon  } \right\|_2^2 } + \frac{1}
{2\pi^2r^2l}\int_T \upsilon dV\right]},
\end{eqnarray*}
and the first part of the lemma is proved.\\

$\bf{2.}$ Let $\upsilon\in C_G^\infty (T)$, with $\int_ {\partial
T} \upsilon dS =0$, $ \phi  = \upsilon  \circ \xi ^{ - 1}$ and $n$
the outward unit normal.\\
 By Stoke's theorem we have
\begin{eqnarray}\label{eq5.28}
  \int_{\partial T} {e^\upsilon  dS} &=&  2\pi r^2 \int_{\partial D}
  {e^\phi  \left( {l + rt} \right)d\sigma _D }   \nonumber \\
& \leqslant & 2\pi r^2 \left( {l + r}
   \right)\int_{\partial D} {e^\phi  d\sigma _D }   \nonumber \\
& = & 2\pi r^2 \left( {l + r} \right)
  \int_D {div(e^\phi  n)dtds}   \nonumber \\
& =& 2\pi r^2 \left( {l + r} \right)
  \int_D {\left[ {divn + n\left( \phi  \right)} \right]e^\phi  dtds}  \nonumber \\
& \leqslant & 2\pi r^2 \left( {l + r}
  \right)\left[ {C_0 \int_D {e^\phi  dtds + \int_D {\left|
  {\nabla \phi } \right|e^\phi  dtds} } } \right],
\end{eqnarray}
where $C_0  = \sup _D \left( {\left| {divn} \right|} \right)$.\\
By (\ref{eq5.28}), and  because of Theorem $4$ of \cite {Che1}
arises
\begin{equation}\label{eq5.29}
\int_{\partial T} {e^\upsilon  dS \leqslant } 2\pi r^2 \left( {l +
r} \right)\left[ {C_0 \tilde C\exp \left( {\mu \left\| {\nabla
\phi } \right\|_2^2 } \right) + \int_D {\left| {\nabla \phi }
\right|e^\phi  dtds} } \right].
\end{equation}
By H\"older's inequality and by Theorem $3$ of \cite {Che3} we
obtain
\begin{equation}\label{eq5.30}
\int_D {\left| {\nabla \phi } \right|e^\phi  dtds} \leqslant
\left\| {\nabla \phi } \right\|_2 \left( {\int_D {e^{2\phi } dtds}
} \right)^{1/2}  \leqslant \tilde C\left\| {\nabla \phi }
\right\|_2 \exp \left( {2\tilde{\mu} \left\| {\nabla \phi }
\right\|_2^2 } \right),
\end{equation}
where $\tilde{\mu}$ a is constant greatest than $1/8\pi$.\\
From the elementary inequality $t \leqslant C_1 \exp \left(
 {\varepsilon_0t^2 } \right),t\geq 0$, $\varepsilon_0
>0$ and $C_1$ a constant with arbitrary $\varepsilon_0 >0$ and $t =
\left\| {\nabla \phi } \right\|_2 $ we obtain
\begin{equation}\label{eq5.31}
 \left\| {\nabla \phi } \right\|_2
\leqslant C_1 \exp \left( {\varepsilon_0 \left\| {\nabla \phi }
\right\|_2^2 } \right).
\end{equation}
Combining inequalities (\ref{eq5.29}), (\ref{eq5.30}) and
(\ref{eq5.31})we obtain
\begin{eqnarray*}
\int_{\partial T} {e^\upsilon  dS }&\leqslant&  2\pi r^2 \left( {l
+ r} \right)\left[ {C_0\tilde C\exp \left( {\tilde{\mu}\left\|
{\nabla \phi } \right\|_2^2 } \right) + \tilde C\left\| {\nabla
\phi }
 \right\|_2 \exp \left( {2\tilde{\mu}\left\| {\nabla \phi } \right\|_2^2 }
  \right)} \right] \hfill \\
 &\leqslant& 2\pi r^2 \left( {l + r}
\right)\\&& \times\left[ {C_0 \tilde C\exp \left( {\tilde{\mu}
\left\| {\nabla \phi }
 \right\|_2^2 } \right)+ \tilde CC_1 \exp \left( {\varepsilon_0 \left\|
 {\nabla \phi } \right\|_2^2 } \right)\exp \left( {2\tilde{\mu}\left\|
 {\nabla \phi }
 \right\|_2^2 } \right)} \right] \hfill \\
 &\leqslant& 2\pi r^2 \left( {l + r}
\right)\tilde C\left( {C_0+ C_1 } \right)\exp \left[ {\left(
{2\tilde{\mu}  +   \varepsilon_0 } \right)\left\| {\nabla \phi }
\right\|_2^2 } \right].
\end{eqnarray*}
Since
$$
 \left\| {\nabla \phi } \right\|_2^2 = \int_D {\left|
{\nabla \phi } \right|^2 dtds}  \leq \frac{1} {{L}}\int_T {\left|
{\nabla \upsilon } \right|} ^2 dV,
$$
the latter inequality becomes
\begin{eqnarray*}
  \int_{\partial T} {e^\upsilon  dS }&\leqslant & 2\pi r^2 \left( {l + r}
   \right)\tilde C\left( {C_0  + C_1 } \right)\exp \left( {\frac{{2\tilde{\mu } + \varepsilon_0 }}
{{L}}\int_T {\left| {\nabla \upsilon } \right|} ^2 dV} \right) \hfill \\
& \leqslant & C \exp \left( {\frac{{2\tilde{\mu } + \varepsilon_0
}} {{L}}\int_T {\left| {\nabla \upsilon } \right|} ^2 dV} \right).
\end{eqnarray*}
Given $\varepsilon>0$, we can choose $\varepsilon_0>0$ such that
$$
\frac{{2\tilde \mu  + \varepsilon _0 }} {{L}} < \frac{2{\tilde \mu
}} {{L}} + \varepsilon =\mu+\varepsilon,
$$
and the last inequality yields
\begin{equation}\label{eq5.32}
\int_{\partial T} {e^\upsilon  dS} \, \leqslant C \exp \left(
{\left( {\mu  + \varepsilon } \right)\int_T {\left| {\nabla
\upsilon } \right|} ^2 dV} \right).
\end{equation}
Rewriting (\ref{eq5.32}) with $ \tilde \upsilon  = \upsilon  -
\frac{1} {4\pi^2r^2l}\int_{\partial T} {\upsilon dS} $ yields
the second inequality of the lemma.\mbox{ }\hfill $\Box$\\

\end{document}